\documentclass[11pt,reqno]{amsart}
\usepackage{amsfonts,amssymb,amsmath,color}

\DeclareMathOperator{\bmin}{\mathbf{min}}
\DeclareMathOperator{\bmax}{\mathbf{max}}
\DeclareMathOperator{\NC}{\mathtt{NC}}
\DeclareMathOperator{\acyclic}{\mathrm{acyclic}}
\DeclareMathOperator{\KG}{\mathtt{KG}}
\DeclareMathOperator{\vvert}{\mathrm{vert}}
\DeclareMathOperator{\conv}{\mathrm{conv}}
\DeclareMathOperator{\precdot}{\prec\!\!\cdot}

\newtheorem{theorem}{Theorem}
\newtheorem{definition}[theorem]{Definition}
\newtheorem{corollary}[theorem]{Corollary}
\newtheorem{lemma}[theorem]{Lemma}
\newtheorem{proposition}[theorem]{Proposition}

\newtheorem{remark}[theorem]{Remark}
\newtheorem{example}[theorem]{Example}
\newtheorem{algorithm}[theorem]{Algorithm}

\numberwithin{equation}{section} \numberwithin{theorem}{section}

\newenvironment{sketch}{\vspace*{2eX}\noindent
{\em Sketch of proof.} }{\hfill$\square$ \\[2eX]}

\begin{document}

\title{Pattern Recognition on Oriented Matroids:
The Existence of a Tope Committee}

\author{Andrey O. Matveev}
\address{Data-Center Co., RU-620034, Ekaterinburg,
P.O.~Box~5, Russian~Federation} \email{aomatveev@\{dc.ru,
hotmail.com\}}


\begin{abstract}
Oriented matroids can serve as a tool of modeling of collective decision-making processes in contradictory
problems of pattern recognition. We present a generalization of the committee techniques of pattern
recognition to oriented matroids. A tope committee for an oriented matroid is a subset of its maximal
covectors such that every positive halfspace contains more than half of the covectors from this subset. For a large subfamily of oriented matroids their committee structure is quite rich; for example, any maximal chains in their tope posets provide one with information sufficient to construct a committee.
\end{abstract}

\maketitle

\pagestyle{myheadings}

\markboth{PATTERN RECOGNITION ON ORIENTED MATROIDS}{A.O.~MATVEEV}

\thispagestyle{empty}

\tableofcontents

\newpage

\section{Introduction}

In this paper we present a description of the pattern recognition
problem in the language of oriented matroids, and show that for
oriented matroids from a large family there exist subsets of their
maximal covectors that can serve as building blocks of collective
decision-making rules.

Book~\cite{DHS} is a standard text on pattern recognition.

In supervised learning to recognize, by means of synthesis of decision rules that make use of decision
surfaces in the feature space, a teacher provides a~class label for every pattern in the training set. The
training patterns from the same class compose a training sample. The union of the training samples is called
the training set. The effect of any decision rule is to divide the~feature space into decision regions each
of which contains training patterns of at most one class. A classifier, that is a pattern recognition system,
relates a new unclassified pattern to a certain class partially presented by a training sample, on the basis
of the inclusion of the new pattern into a~decision region; the decision rule also can force the classifier
to leave the new pattern unclassified.

If decision surfaces in the feature space $\mathbb{R}^n$ are hyperplanes, then decision regions are open
convex polyhedra.

The theory of hyperplane arrangements is an area of an active
cross-disciplinary study, see, e.g.,~\cite{OT,St-HA}.

If a training set in $\mathbb{R}^n$ is composed of two samples
which cannot be separated by one decision hyperplane then a
classifier (called in the case of two classes a dichotomizer)
operates with decision rules relying on infeasible linear
inequality systems.

Various properties of infeasible linear constraint systems have
been studied in depth, see, e.g.,
\cite{A-Thesis,Ch,E1,E2,EM,EMA,EMSKh,M,Pa-Thesis,P-Thesis}.

A useful generalization of the notion of solution of a linear inequality system to the infeasible case is the
notion of majority committee.

A committee for an infeasible system of strict linear inequalities
over $\mathbb{R}^n$ is a finite subset of elements of
$\mathbb{R}^n$ such that for every inequality more than half of
the elements are its solutions.

Such committees were apparently first introduced in seminal
notes~\cite{AK-IEEE,AK-AMS}, where their application to the
pattern recognition problem was discussed. Various committee
constructions for a variety of contradictory problems have later
been invented and explored in detail; some of the surveys in this
subject are
works~\cite{EMSKh,Kh-DSciThesis,KhMR,M,MKBKS,MKh1,MKh2}.

Infeasible systems of strict linear inequalities over
$\mathbb{R}^n$, real hyperplane arrangements for which the
intersections of all positive halfspaces of $\mathbb{R}^n$ are
empty, and realizable oriented matroids which are not acyclic, are
closely connected mathematical objects that allow one to model
mechanisms of collective decision-making in pattern recognition
problems posed in terms of infeasible systems of linear
constraints.

Oriented matroids are defined by various equivalent axiom systems,
and they can be thought of as a combinatorial abstraction of point
configurations over the reals, of real hyperplane arrangements, of
convex polytopes, and of directed graphs.

Oriented matroids are reviewed and studied in detail, e.g.,
in~\cite{BK,BLSWZ,JGB,R-GZ,S,Z}.

In the present paper, a direct generalization of the notion of committee for a linear inequality system, in
terms of maximal covectors of an oriented matroid, is presented: A tope committee $\mathcal{K}^{\ast}$ for an
oriented matroid is a subset of its maximal covectors such that for every element of the ground set the
corresponding positive halfspace of the oriented matroid contains more than half of the covectors from
$\mathcal{K}^{\ast}$.

One of the approaches to the study of committees consists in structural and combinatorial analysis of the
family of maximal feasible subsystems of constraints, and in investigation of the properties of graphs which
are naturally associated with those subsystems~\cite{GNT}; such graphs are an example of constructions dual
to the Kneser graphs of set systems considered \mbox{in~\cite[\S{}3.3]{MBZ},\cite{MZ}.} The properties of the
graphs associated with the maximal feasible subsystems, such as connectedness and the existence of an odd
cycle, are important for graph-theoretic algorithms of synthesis of committees.

The present paper constitutes a review of central ideas from
works~\cite{AK-IEEE} and~\cite{GNT}, formulated in the language of
oriented matroids.

In Section~\ref{section:1}, some terminology of the theory of oriented matroids used in the paper is
recalled. In Section~\ref{section:2}, we review the setting of the pattern recognition problem, give the
definitions of committees of maximal (co)vectors, and discuss committee decision rules. In central
Section~\ref{section:3}, it is shown that every oriented matroid without loops and antiparallel elements, has
a tope committee of cardinality less than or equal to the cardinality of its ground set. The argument is
based on analysis of consecutive reorientations of an initial acyclic oriented matroid.
Section~\ref{section:4} is devoted to graphs which are naturally associated with the families of topes with
inclusion-maximal positive parts. The sets of vertices of the odd cycles in such graphs are committees for
the corresponding oriented matroids. We apply the graph-theoretic approach from~\cite{GNT} to a
generalization of the basic construction of centrally-symmetric cycle of adjacent regions in a hyperplane
arrangement from~\cite{AK-IEEE}; and vice versa, we use `symmetric cycles' in the tope graph of an oriented
matroid, inspired by the above-mentioned cycles of regions from~\cite{AK-IEEE}, to prove several generalized
graph-theoretic results from~\cite{GNT}. Section~\ref{section:5} mentions the link between the committees for
an oriented matroid and blocker constructions in the Boolean lattice of subsets of the tope set.

\section{Preliminaries}
\label{section:1}

All oriented matroids considered in the paper are of rank that is greater than or equal to $2$. We use quite
nonstandard definitions of simple oriented matroids and of graph homomorphisms: An oriented matroid  is {\em
simple\/} if it has no loops, parallel or {\sl antiparallel\/} elements. A {\em homomorphism\/} of a graph to
a graph is a mapping from the vertex set of the first graph to that of the second graph, such that either the
image of any edge is an edge, {\sl or\/} the images of the endvertices coincide.

$E_m$ denotes the set $[1,m]:=\{1,2,\ldots,m\}$. $\mathrm{T}^{(+)}$ denotes the sign vector $(++\ldots+)$
whose components are all $+$; $\mathrm{T}^{(-)}:=-\mathrm{T}^{(+)}=(--\ldots-)$.

See~\cite[Chapter~3]{St-EC},\cite[\S{}4.1]{BLSWZ} on posets. If
$X$ is a subset of a poset $\mathfrak{P}$, then $\bmin X $ denotes
the set of all minimal elements from $X$. $\mathfrak{I}(X)$ and
$\mathfrak{F}(X)$ denote the order ideal and filter in
$\mathfrak{P}$ generated by $X$, respectively. If $\mathfrak{P}$
is graded then $\mathfrak{P}^{(k)}$ is the set of all its elements
of poset rank $k$.

For a set family $\mathcal{F}:=\{F_i:\ i\in[1,m]\}$, $\bmin\mathcal{F}$ and $\bmax\mathcal{F}$ denote the
subfamilies of all inclusion-minimal and of all inclusion-maximal sets in $\mathcal{F}$, respectively. The
{\em nerve\/} of $\mathcal{F}$ is an abstract simplicial complex on the vertex set $\mathcal{F}$; a subset
$K\subseteq[1,m]$ is a face of the nerve iff $|\bigcap_{k\in K}F_k|>0$, see,
e.g.,~\cite[\S{}10]{B-TM},\cite[\S{}8.5]{P}.

See~\cite{Diestel} on graphs. Throughout the paper, graphs are undirected; they have no loops and multiple
edges. For a graph $\mathbf{G}$, its sets of vertices and of edges are denoted by $\mathfrak{V}(\mathbf{G})$
and $\mathfrak{E}(\mathbf{G})$, respectively. {\em Cycles\/} are regular subgraphs of valency $2$; all
vertices of {\em paths\/} in graphs are distinct.

The {\em neighborhood complex\/} $\NC(\mathbf{G})$ of a graph
$\mathbf{G}$, defined in~\cite{L}, is an abstract simplicial
complex on the vertex set $\mathfrak{V}(\mathbf{G})$; a subset
$N:=\{n_1,\ldots,n_k\}$ $\subset\mathfrak{V}(\mathbf{G})$ is a
face of the complex iff there is $v\in\mathfrak{V}(\mathbf{G})$
such that $\{n_1,v\},\ldots,$
$\{n_k,v\}\in\mathfrak{E}(\mathbf{G})$.

Recall that for graphs $\mathbf{G}'$ and $\mathbf{G}''$, a {\em
homomorphism\/} of $\mathbf{G}'$ to $\mathbf{G}''$, written as
$h:\mathbf{G}'\to \mathbf{G}''$, is a mapping
$h:\mathfrak{V}(\mathbf{G}')\to\mathfrak{V}(\mathbf{G}'')$ such
that $\{\mathrm{u},\mathrm{v}\}\in\mathfrak{E}(\mathbf{G}')$
implies
$\{h(\mathrm{u}),h(\mathrm{v})\}\in\mathfrak{E}(\mathbf{G}'')$ or
$h(\mathrm{u})=h(\mathrm{v})$.

If $\mathcal{F}$ is a set family then the {\em Kneser graph\/} $\KG(\mathcal{F})$ of $\mathcal{F}$,
considered in~\cite[\S{}3.3]{MBZ},\cite{MZ}, is the graph with
$\mathfrak{V}\bigl(\KG(\mathcal{F})\bigr):=\mathcal{F}$; if $F',F''\in\mathcal{F}$ then $\{F',F''\}$
\mbox{$\in\mathfrak{E}\bigl(\KG(\mathcal{F})\bigr)$} iff $|F'\cap F''|=0$.

We borrow almost all terminology concerning oriented matroids
from~\cite[Chapters~3, 4, 7]{BLSWZ}:

Let $E$ be a finite set, $\{-,0,+\}$ the set of signs, and $\{-,0,+\}^E$ the set of {\em sign vectors}. The
{\em support\/} of a sign vector $X\in\{-,0,+\}^E$ is $\underline{X}:=\{e$ $\in E:\ X(e)\neq 0\}$; here
$X(e)$ denotes the $e$th component of $X$. $X^-:=\{e$ $\in E:\ X(e)=-\}$ denotes the set of {\em negative
elements\/} of $X$; $X^+:=\{e\in E:\ X(e)=+\}$ is the set of {\em positive elements\/} of $X$. Thus
$\underline{X}:=X^-\cup X^+$. $X^-$ and $X^+$ are also called the {\em negative\/} and {\em positive parts\/}
of $X$, respectively. An inclusion $e\in X$ means $e\in\underline{X}$. The {\em zero\/} sign vector
$(00\ldots0)$, with the empty support, is denoted by $\pmb{0}$. The {\em zero set\/} $\mathbf{z}(X)$ of a
sign vector $X$ is the set $\{e\in E:\ X(e)=0\}$.

If $\mathcal{P}$ is a set of sign vectors then $\bmax^+(\mathcal{P}):=\bigl\{P\in\mathcal{P}:\
P^+\in\bmax\{R^+:\ R\in\mathcal{P}\}\bigr\}$; similarly, $\bmin^+(\mathcal{P}):=\bigl\{P\in\mathcal{P}:\
P^+\in\bmin\{R^+:\ R\in\mathcal{P}\}\bigr\}$ $=\{-P: P\in\bmax^+(\mathcal{P})\}$.

If $A\subseteq E$ then the sign vector ${}_{-A}X$ is defined by
\begin{equation*}
({}_{-A}X)(e):=\begin{cases}+,&\text{if $e\in A$ and $X(e)=-$}\
,\\ -,&\text{if $e\in A$ and $X(e)=+$}\ ,\\
X(e),&\text{otherwise}\end{cases}
\end{equation*}
(if $e\in E$ then we write ${}_{-e}X$ instead of ${}_{-\{e\}}X$).
In particular, the {\em opposite\/} of $X$ is $-X:={}_{-E}X$, that
is,
\begin{equation*}
(-X)(e):=\begin{cases}+,&\text{if $X(e)=-$}\ ,\\ -,&\text{if
$X(e)=+$}\ ,\\ 0,&\text{if $X(e)=0$}\ ,\end{cases}
\end{equation*}
for all $e\in E$. If $\mathcal{F}\subseteq\{-,0,+\}^E$ and $A\subseteq E$, then
${}_{-A}\mathcal{F}:=\{{}_{-A}X:\ X\in\mathcal{F}\}$; in particular, $-\mathcal{F}:={}_{-E}\mathcal{F}=\{-X:\
X\in\mathcal{F}\}$.

If $X\in\{-,0,+\}^E$ then the sign vector $X$ is called {\em
nonpositive\/} (resp., {\em negative}) if $X(e)\in\{-,0\}$ (resp.,
$X(e)=-$), for all $e\in E$. Similarly, $X$ is {\em nonnegative\/}
(resp., {\em positive}) if $-X$ is {\em nonpositive\/} (resp.,
{\em negative}).

The {\em composition\/} of two sign vectors $X$ and $Y$ is the
sign vector $X\circ Y$ defined by
\begin{equation*}
(X\circ Y)(e):=\begin{cases}X(e),&\text{if $X(e)\neq 0$}\ ,\\
Y(e),&\text{otherwise}\ .\end{cases}
\end{equation*}

The {\em separation set\/} of $X$ and $Y$ is $\mathbf{S}(X,Y):=\{e\in E:\ X(e)=-Y(e)\neq 0\}$. If
$|\mathbf{S}(X,Y)|=0$ then one says that the sign vectors $X$ and $Y$ are {\em conformal}; in this case
$X\circ Y=Y\circ X$. If sign vectors $X_1,X_2,\ldots,X_k\in\{-,0,+\}^E$ are pairwise conformal then
$\underset{i\in[1,k]}\bigcirc X_i$ is a short notation for the {\em conformal composition\/} $X_1\circ
X_2\circ\cdots\circ X_k$.

The partial order on the set $\{-,0,+\}$ is defined by the
relations $0<-$ and $0<+$; the signs $-$ and $+$ are incomparable.
This induces the product partial order on $\{-,0,+\}^E$, in which
sign vectors are compared componentwise. Thus $X\leq Y$ iff
$X(e)\in\{0,Y(e)\}$ for all $e\in E$.

Oriented matroids are defined by several equivalent axiom systems.

Let $E$ be a finite set. If $\mathcal{C}\subseteq\{-,0,+\}^E$,
then $\mathcal{C}$ by definition is the set of {\em circuits\/} of
an oriented matroid on $E$ iff it satisfies the following {\em
Circuit Axioms}~\cite[Definition~3.2.1]{BLSWZ}:
\begin{itemize}
\item[\rm(C0)] $\pmb{0}\not\in\mathcal{C}$;
\item[\rm(C1)] $X\in\mathcal{C}$ implies $-X\in\mathcal{C}$;
\item[\rm(C2)] $X,Y\in\mathcal{C}$ and
$\underline{X}\subseteq\underline{Y}$ imply $X=Y$ or $X=-Y$;
\item[\rm(C3)] if $X,Y\in\mathcal{C}$, $X\neq-Y$, and
$e\in X^+\cap Y^-$, then there is $Z\in\mathcal{C}$ such that
$Z^-\subseteq (X^-\cup Y^-)-\{e\}$ and $Z^+\subseteq (X^+\cup
Y^+)-\{e\}$.
\end{itemize}
An oriented matroid on $E$, with set of circuits $\mathcal{C}$, is
denoted by $(E,\mathcal{C})$.

The circuit supports $\underline{\mathcal{C}}:=\{\underline{C}:\ C\in\mathcal{C}\}$ in an oriented matroid
$\mathcal{M}$ $:=(E,\mathcal{C})$ constitute the circuits of the {\em underlying matroid of\/} $\mathcal{M}$,
denoted by $\underline{\mathcal{M}}$. The {\em rank\/} of $\mathcal{M}$ by definition is the rank of
$\underline{\mathcal{M}}$.

A {\em vector\/} of an oriented matroid is any composition of its
circuits. An oriented matroid on $E$, given by set of its vectors
$\mathcal{V}$, is denoted by $(E,\mathcal{V})$. A {\em maximal
vector\/} of an oriented matroid is a vector whose support is
maximal with respect to inclusion. An oriented matroid on a set
$E$, with set of maximal vectors $\mathcal{W}$, is denoted by
$(E,\mathcal{W})$.

If $\mathcal{\mathcal{L}}\subseteq\{-,0,+\}^E$, then the pair
$(E,\mathcal{L})$ is an oriented matroid on $E$, with the set of
{\em covectors\/} $\mathcal{L}$, iff $\mathcal{L}$ satisfies the
following {\em Covector Axioms}~\cite[Proposition~4.1.1]{BLSWZ}:
\begin{itemize}
\item[\rm(L0)] $\pmb{0}\in\mathcal{L}$;
\item[\rm(L1)] $X\in\mathcal{L}$ implies $-X\in\mathcal{L}$;
\item[\rm(L2)] $X,Y\in\mathcal{L}$ implies $X\circ
Y\in\mathcal{L}$;
\item[\rm(L3)] if $X,Y\in\mathcal{L}$ and $e\in\mathbf{S}(X,Y)$
then there exists $Z\in\mathcal{L}$ such that $Z(e)=0$ and
$Z(f)=(X\circ Y)(f)=(Y\circ X)(f)$ for all
$f\not\in\mathbf{S}(X,Y)$.
\end{itemize}

For an oriented matroid $(E,\mathcal{L})$ of rank $r$, the poset
$\widehat{\mathcal{L}}:=\mathcal{L}\dot\cup\{\hat{1}\}$, with a top element $\hat{1}$ adjoined, is a graded
lattice (the so-called {\em `big' face lattice}) of length $r+1$, see~\cite[Theorem~4.1.14]{BLSWZ};
$\hat{0}=\mathbf{0}$ is the bottom element of $\widehat{\mathcal{L}}$.

A {\em maximal covector\/} (a {\em tope}) of an oriented matroid
is a covector whose support is maximal with respect to inclusion.
An oriented matroid $\mathcal{M}$ on a set $E$, with set of topes
$\mathcal{T}$, is denoted by $(E,\mathcal{T})$.

The set $\mathcal{C}^{\ast}$ of non-zero covectors of an oriented matroid $\mathcal{M}$, with
inclusion-minimal supports, is the set of {\em cocircuits\/} of $\mathcal{M}$. An oriented matroid on $E$,
with that set of cocircuits, is denoted by $(E,\mathcal{C}^{\ast})$.

For every $e\in E$ the corresponding {\em positive halfspace\/} is the subset of topes
$\mathcal{T}_e^+:=\{T\in\mathcal{T}:\ T(e)=+\}$; the {\em negative halfspace\/} $\mathcal{T}_e^-$ is the
subset $-\mathcal{T}_e^+$.

The set of vertices of the {\em tope graph\/} $\mathcal{T}(\mathcal{L})$ is the set of topes; two topes are
connected by an edge if the topes are {\em adjacent}, that is, if they cover the same element (a {\em
subtope}) of poset rank $r-1$ in $\widehat{\mathcal{L}}$, where $r$ is the rank of the oriented matroid.

If $B\in\mathcal{T}$ then the {\em tope poset\/}
$\mathcal{T}(\mathcal{L},B)$, based at $B$, is defined by the
partial order on the set of topes: $T'\preceq T''$ iff
$\mathbf{S}(B,T')\subseteq\mathbf{S}(B,T'')$.

A subset of topes $\mathcal{Q}\subseteq\mathcal{T}$ is {\em
$\mathrm{T}$-convex\/} if the following implication holds:
\begin{equation*}
T',T''\in\mathcal{Q}\ ,\ \ T\in\mathcal{T}\ ,\ \
|\mathbf{S}(T',T'')|=|\mathbf{S}(T',T)|+|\mathbf{S}(T,T'')|\ \
\Longrightarrow\ \ T\in\mathcal{Q} ,
\end{equation*}
that is, if $\mathcal{Q}$ contains every shortest path in the
graph $\mathcal{T}(\mathcal{L})$ between any two of its members.
The {\em $\mathrm{T}$-convex hull\/}
$\conv_{\mathrm{T}}(\mathcal{Q})$ of
$\mathcal{Q}\subseteq\mathcal{T}$ is the intersection of all
halfspaces that contain $\mathcal{Q}$.

All topes $T\in\mathcal{T}$ have the same support and the same
zero set $E_{\circ}:=\mathbf{z}(T)$. The elements in $E_{\circ}$
are called the {\em loops\/} of $\mathcal{M}$. Thus $e\in E$ is a
loop of $\mathcal{M}:=(E,\mathcal{C})$ iff there is a circuit
$(0\ldots 0 \underset{\substack{\uparrow\\ e}}{+} 0 \ldots
0)\in\mathcal{C}$.

If $e\not\in C$ for every circuit $C\in\mathcal{C}$ of an oriented
matroid $\mathcal{M}:=(E,\mathcal{C})$, then $e$ is called a {\em
coloop\/} of $\mathcal{M}$.

Elements $e,f\in E$, $e\neq f$, are called {\em parallel\/} if
$X(e)=X(f)$ for all $X\in\mathcal{L}$; they are called {\em
antiparallel}, if $X(e)=-X(f)$ for all $X\in\mathcal{L}$.

As mentioned above, a {\em simple
oriented matroid\/} means, throughout the paper, an oriented matroid without loops,
parallel or antiparallel elements.

The {\em restriction\/} of a sign vector $X\in\{-,0,+\}^E$ to a
subset $A\subseteq E$ is the sign vector $X|_A\in\{-,0,+\}^A$
defined by $(X|_A)(e):=X(e)$ for all $e\in A$.

For an oriented matroid $\mathcal{M}:=(E,\mathcal{L})$, the oriented matroid $(E-A,\mathcal{L}\backslash A)$
on $E-A$, given by its set of covectors $\mathcal{L}\backslash A:=\{X|_{E-A}:\ X\in\mathcal{L}\}$
\mbox{$\subseteq\{-,0,+\}^{E-A}$} is called the {\em deletion\/} $\mathcal{M}\backslash A$ or the {\em
restriction\/} $\mathcal{M}|_{E-A}$. The oriented matroid $(E,{}_{-A}\mathcal{L})$ on $E$, given by its set
of covectors ${}_{-A}\mathcal{L}\subseteq\{-,0,+\}^E$ is called the {\em reorientation\/}
${}_{-A}\mathcal{M}$; see~\cite[Lemma~4.1.8]{BLSWZ}.

An oriented matroid $\mathcal{M}:=(E,\mathcal{C})=(E,\mathcal{T})$
is {\em acyclic} if there is no nonnegative circuit in the set
$\mathcal{C}$, or equivalently if there exists the nonnegative
tope in $\mathcal{T}$. A subset $A\subseteq E$ is called {\em
acyclic\/} if the restriction $\mathcal{M}|_{A}$ is acyclic. The
oriented matroid $\mathcal{M}$ is {\em totally cyclic\/} if for
each element $e\in E$ there exists a nonnegative circuit
$C\in\mathcal{C}$ such that $e\in C$. Recall that `most' oriented
matroids are neither acyclic nor totally
cyclic~\cite[\S{}6.3.1]{R-GZ}.

The circuits, vectors, and maximal vectors of an oriented matroid $\mathcal{M}$ are the cocircuits,
covectors, and topes, respectively, of the oriented matroid $\mathcal{M}^{\ast}$, the {\em dual\/} (or {\em
orthogonal}\/) of $\mathcal{M}$. The loops of $\mathcal{M}$ are the coloops of $\mathcal{M}^{\ast}$;
$\mathcal{M}$ is acyclic iff $\mathcal{M}^{\ast}$ is totally cyclic, see~\cite[Proposition~3.4.8]{BLSWZ}.

For an oriented matroid $\mathcal{M}:=(E,\mathcal{C}^{\ast})$, a {\em single element extension\/}
$\widetilde{\mathcal{M}}$ $:=(\widetilde{E},\widetilde{\mathcal{C}^{\ast}})$ of $\mathcal{M}$ is an oriented
matroid on a set $\widetilde{E}$ such that $\widetilde{E}=E\dot\cup\{g\}$, with set of cocircuits
$\widetilde{\mathcal{C}^{\ast}}$. If $g$ is not a coloop of $\widetilde{\mathcal{M}}$, then
$\widetilde{\mathcal{M}}$ is called a {\em nontrivial\/} extension of $\mathcal{M}$.

The set $\widetilde{\mathcal{C}^{\ast}}$ of cocircuits of an
extension $\widetilde{\mathcal{M}}$ is described in the following
way~\cite{LasV},~\cite[Proposition~7.1.4]{BLSWZ}:
\begin{itemize}
\item[\rm(i)] Let $\widetilde{\mathcal{M}}$ be a nontrivial single element
extension of $\mathcal{M}:=(E,\mathcal{C}^{\ast})$ $=(E,\mathcal{L})$. Then for every cocircuit
$Y\in\mathcal{C}^{\ast}$ there is a unique way to extend $Y$ to a cocircuit of $\widetilde{\mathcal{M}}$:
there is a unique function $\sigma:\mathcal{C}^{\ast}$ $\to\{-,0,+\}$, called the {\em localization}, such
that $\bigl\{\bigl(Y,\sigma(Y)\bigr):\ Y$
\mbox{$\in\mathcal{C}^{\ast}\bigr\}\subseteq\widetilde{\mathcal{C}^{\ast}}$}, that is,
$\bigl(Y,\sigma(Y)\bigr)$ is a cocircuit of $\mathcal{M}$ for every cocircuit $Y$ of $\mathcal{M}$.
Furthermore, this $\sigma$ satisfies $\sigma(-Y)=-\sigma(Y)$ for all $Y$ $\in\mathcal{C}^{\ast}$.
\item[\rm(ii)] $\widetilde{\mathcal{M}}$ is uniquely determined by
$\sigma$, with
\begin{multline*}
\widetilde{\mathcal{C}^{\ast}}=\bigl\{\ \bigl(Y,\sigma(Y)\bigr):\
Y\in\mathcal{C}^{\ast}\ \bigr\}\\ \dot\cup\ \ \Bigl\{\ (Y'\circ
Y'',0):\ Y',Y''\in\mathcal{C}^{\ast},\ \sigma(Y')=-\sigma(Y'')\neq
0,\\ |\mathbf{S}(Y',Y'')|=0,\ \rho(Y'\circ Y'')=2\ \Bigr\}\ ,
\end{multline*}
where $\rho$ denotes the poset rank function on
$\widehat{\mathcal{L}}$.
\end{itemize}

\section{Pattern Recognition on Oriented Matroids}
\label{section:2}
\subsection{The Two-Class Pattern Recognition Problem}

A {\em training set\/} is a simple oriented matroid $\mathcal{S}$
on a ground set $E$, together with a mapping $\lambda:E\to\{-,+\}$
such that the {\em training samples\/} $\lambda^{-1}(-)$ and
$\lambda^{-1}(+)$ are nonempty. The elements of $E$ are called the
{\em training patterns}.

{\em Classes\/} $\mathbf{A}$ and $\mathbf{B}$ are disjoint sets such that
$\mathbf{A}\supseteq\lambda^{-1}(-)$ and $\mathbf{B}\supseteq\lambda^{-1}(+)$. Thus, an element $e\in E$ {\sl
a priori\/} belongs to the class $\mathbf{A}$ iff $\lambda(e)=-$; it {\sl a priori\/} belongs to the class
$\mathbf{B}$ iff $\lambda(e)=+$.

Let $\widetilde{\mathcal{S}}$ denote the nontrivial single element
extension of $\mathcal{S}$ by a new unclassified pattern $g$ which
is not a loop, and which is parallel or antiparallel to neither of
the elements of $E$. A {\em decision rule\/} is any mapping
\begin{equation*}
\mathfrak{r}:\ E\dot\cup\{g\}\to\{-,0,+\}
\end{equation*}
such that
\begin{equation*}
\mathfrak{r}:\ e\mapsto \lambda(e)\ ,\ \ \ e\in E\ .
\end{equation*}
If $f\in E\dot\cup\{g\}$ and $\mathfrak{r}(f)=-$, then a
dichotomizer using the rule $\mathfrak{r}$ relates the pattern $f$
to the class $\mathbf{A}$. If $\mathfrak{r}(f)=+$, then $f$ is
classified by the dichotomizer as a pattern from $\mathbf{B}$. If
$\mathfrak{r}(f)=0$ then the dichotomizer leaves the pattern $f$
unclassified.

\subsection{Committees of Maximal (Co)Vectors for Oriented
Matroids}

\begin{definition}
\label{def:2} Let $p$ be a rational number such that $0\leq p<1$.
\begin{itemize}
\item
Given an oriented matroid $\mathcal{M}:=(E,\mathcal{T})$, a subset
$\mathcal{K}^{\ast}\subset\mathcal{T}$ is a {\em tope
$p$-committee} {\rm(}a {\em $p$-committee of maximal
covectors}{\rm)} {\em for $\mathcal{M}$} if for every $e\in E$ it
holds
\begin{equation*}
|\{K\in\mathcal{K}^{\ast}:\ K(e)=+\}|>p|\mathcal{K}^{\ast}|\ .
\end{equation*}

A tope $\tfrac{1}{2}$-committee for $\mathcal{M}$ is called a {\em
tope committee for $\mathcal{M}$}.

\item
Given an oriented matroid $\mathcal{M}:=(E,\mathcal{W})$, a subset
$\mathcal{K}\subset\mathcal{W}$ is a {\em $p$-committee of maximal
vectors for $\mathcal{M}$} if for every $e\in E$ it holds
\begin{equation*}
|\{K\in\mathcal{K}:\ K(e)=+\}|>p|\mathcal{K}|\ .
\end{equation*}

A $\tfrac{1}{2}$-committee of maximal vectors for $\mathcal{M}$ is
called a {\em committee of maximal vectors for} $\mathcal{M}$.
\end{itemize}
\end{definition}

Definition~\ref{def:2} implies that
\begin{itemize}
\item
a set $\mathcal{K}^{\ast}$ is a tope $p$-committee for
$\mathcal{M}$ iff $\mathcal{K}^{\ast}$ is a $p$-committee of
maximal vectors for the dual oriented matroid
$\mathcal{M}^{\ast}$;
\item
a subset $\mathcal{K}^{\ast}\subset\mathcal{T}$ is a committee for
$\mathcal{M}$ iff for every $e\in E$ it holds
$|\mathcal{K}^{\ast}|<2|\{K\in\mathcal{K}^{\ast}:\ K(e)=+\}|$;
\item
a subset $\mathcal{K}^{\ast}\subset\mathcal{T}$ is a committee for
$\mathcal{M}$ iff the set $\{-T:\
T\in\mathcal{T}-\mathcal{K}^{\ast}\}$ is;
\item
if $\mathcal{K}^{\ast}$ is a tope committee for $\mathcal{M}$, and
if $\widetilde{\mathcal{M}}$ is a trivial single element extension
of $\mathcal{M}$ by a coloop, then the set $\bigl\{(T,+):\
T\in\mathcal{K}^{\ast}\bigr\}$ is a tope committee for
$\widetilde{\mathcal{M}}$.
\end{itemize}

Given an oriented matroid $\mathcal{M}$, we denote by
$\mathbf{K}^{\ast}(\mathcal{M})$ the family of all tope committees
for $\mathcal{M}$.

\begin{definition}
Let $\mathcal{M}:=(E,\mathcal{T})$ be an oriented matroid, and
$\mathcal{K}^{\ast}$ a tope committee for $\mathcal{M}$.
\begin{itemize}
\item
$\mathcal{K}^{\ast}$ is called  {\em minimal\/} if any proper
subset of the set $\mathcal{K}^{\ast}$ is not a committee for
$\mathcal{M}$.
\item
If $\mathcal{K}^{\ast}$ is minimal, then $\mathcal{K}^{\ast}$ is
called {\em critical\/} if
\begin{equation*}
K\in\mathcal{K}^{\ast}\ ,\ \ T\in\mathcal{T}\ ,\ \ T^+\subsetneqq
K^+\ \ \Longrightarrow\ \ (\mathcal{K}^{\ast}-\{K\})\ \cup\
\{T\}\not\in\mathbf{K}^{\ast}(\mathcal{M})\ .
\end{equation*}
\item
$\mathcal{K}^{\ast}$ is called a {\em minimum committee\/} {\rm(}a
{\em committee of minimal cardinality}{\rm)} if there is no
committee $\mathcal{Q}^{\ast}$ for $\mathcal{M}$ such that
$|\mathcal{Q}^{\ast}|<|\mathcal{K}^{\ast}|$.
\end{itemize}
\end{definition}

If $\mathcal{M}$ is simple and acyclic, then the one-element set
$\{\mathrm{T}^{(+)}\}$ is a critical tope committee for
$\mathcal{M}$.

It follows from the definition that minimal and minimum committees
do not contain opposites.

\begin{example}
Consider the acyclic oriented matroid
$\mathcal{N}^0:=(E_6,\mathcal{T}^0)$ given by the central
hyperplane arrangement of Figure~{\rm\ref{fig:7}}.

\begin{figure}[ht]
\begin{picture}(10,13.7)(0.8,1)
\put(1,4.5){\line(1,2){4}} \put(1,4.5){\line(5,-2){1.65}}
\put(2,2.5){\line(1,2){4}} \put(2,2.5){\line(3,2){1.5}}
\put(3,1){\line(1,5){0.5}} \put(3,1){\line(1,2){4}}
\put(3.5,3.5){\line(1,2){0.7}}

\put(6,8.5){\line(1,2){1.5}}

\put(6.05,2.5){\line(-5,2){1.9}} \put(6.05,2.5){\line(1,2){4}}
\put(2.5,7.5){\line(2,-3){2.15}} \put(2.5,7.5){\line(1,0){2}}
\put(3.5,9.5){\line(2,-3){2.15}} \put(5.65,6.3){\line(4,-1){2}}
\put(4.3,11){\line(2,-3){2.1}} \put(4.65,4.25){\line(6,5){3.9}}
\put(6.4,7.8){\line(1,-6){0.58}} \put(8.5,7.5){\line(-1,0){2.05}}
\put(6,7){\line(2,-3){0.7}} \put(7,4.35){\line(-1,2){0.65}}
\put(7.7,5.85){\line(-1,1){0.5}} \put(4.85,8.15){\line(-1,1){1.3}}
\put(4.3,11){\line(2,-5){0.86}} \put(3.7,5.8){\line(1,1){1.35}}
\put(4.7,4.4){\line(1,4){0.6}} \put(5.45,7.5){\circle * {0.2}}
\put(5.5,7.5){\line(-3,1){0.85}} \put(5.45,7.5){\line(1,-6){0.19}}
\put(5.45,7.5){\line(0,1){1.8}} \put(5.5,7.5){\line(5,2){0.9}}
\put(5,12.5){\line(5,-2){5}} \put(6,10.5){\line(3,2){3}}
\put(9,12.5){\line(-1,-2){0.65}} \put(7,9){\line(1,5){1}}
\put(8,14){\line(-1,-2){1.13}}
\put(9.5,10.87){\vector(-1,-2){0.15}}
\put(7.1,9.5){\vector(4,-1){0.15}}
\put(6.4,10.8){\vector(1,-1){0.15}}
\put(4.8,10.3){\vector(2,1){0.15}} \put(4,8.7){\vector(1,1){0.15}}
\put(3,6.7){\vector(1,1){0.15}}

\put(5.8,10.4){\makebox(0,0)[r]{$\mathbf{1}$}}
\put(7,8.7){\makebox(0,0)[l]{$\mathbf{4}$}}
\put(10.15,10.3){\makebox(0,0)[l]{$\mathbf{3}$}}
\put(4.2,11.2){\makebox(0,0)[r]{$\mathbf{6}$}}
\put(3.45,9.7){\makebox(0,0)[r]{$\mathbf{2}$}}
\put(2.45,7.7){\makebox(0,0)[r]{$\mathbf{5}$}}
\end{picture}
\caption{A central hyperplane arrangement that realizes a simple
acyclic oriented matroid $\mathcal{N}^0:=(E_6,\mathcal{T}^0)$ with
$28$ topes. The positive halfspaces of $\mathbb{R}^3$ are marked
by arrows.} \label{fig:7}
\end{figure}
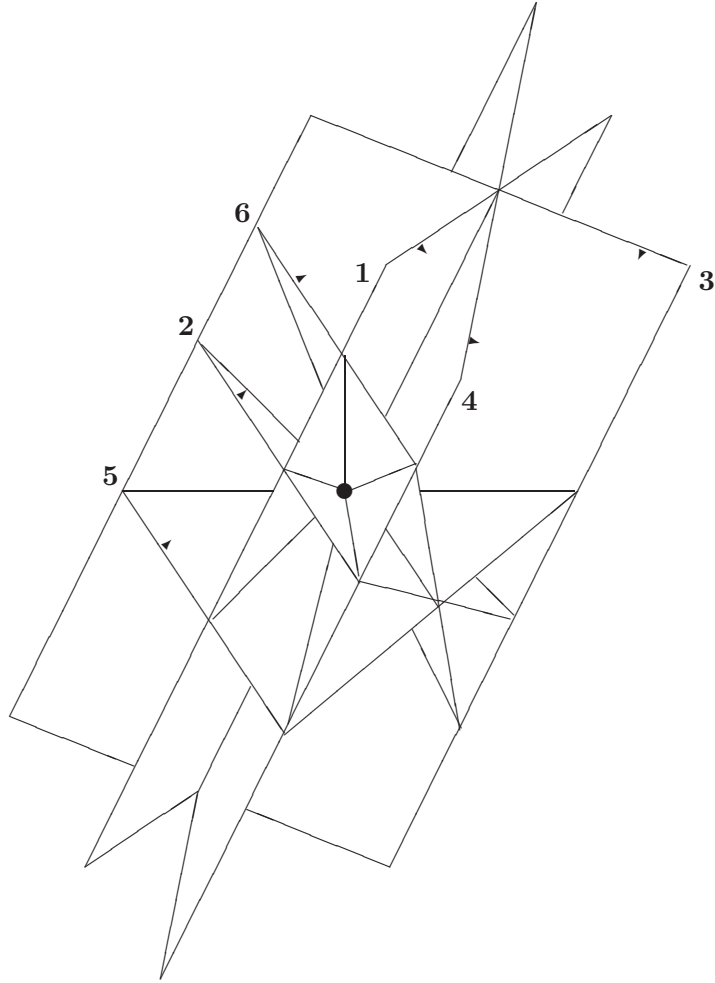

The sets of topes $\mathcal{T}^0$ and $\mathcal{T}^2$ of the oriented matroids $\mathcal{N}^0$ and
$\mathcal{N}^2:={}_{-[1,2]}\mathcal{N}^0=(E_6,\mathcal{T}^2)$, respectively, are as follows:

{\tiny
\begin{equation*}
\mathcal{T}^0:=\quad
\begin{matrix}
\{&+&+&+&+&+&+&\\ &+&+&+&-&+&+&\\ &-&+&+&-&+&+&\\ &-&+&+&-&+&-&\\
&+&+&+&-&+&-&\\ &+&+&+&+&+&-&\\ &+&+&+&+&-&+&\\ &+&-&+&+&-&+&\\
&+&-&+&+&+&-&\\ &+&-&+&-&+&-&\\ &-&-&+&-&+&-&\\ &-&-&+&-&-&-&\\
&+&-&+&-&-&-&\\ &+&-&+&+&-&-&\\ &-&+&-&-&+&+&\\ &-&+&-&+&+&+&\\
&+&+&-&+&+&+&\\ &+&+&-&+&-&+&\\ &-&+&-&+&-&+&\\ &-&+&-&-&-&+&\\
&-&+&-&-&+&-&\\ &-&-&-&-&+&-&\\ &-&-&-&-&-&+&\\ &-&-&-&+&-&+&\\
&+&-&-&+&-&+&\\ &+&-&-&+&-&-&\\ &-&-&-&+&-&-&\\ &-&-&-&-&-&-&\}
\end{matrix}\quad\quad\quad
\mathcal{T}^2=\quad
\begin{matrix}
\{&-&-&+&+&+&+\\ &-&-&+&-&+&+\\ &+&-&+&-&+&+\\ &+&-&+&-&+&-&\\
&-&-&+&-&+&-\\ &-&-&+&+&+&-\\ &-&-&+&+&-&+\\ &-&+&+&+&-&+&\\
&-&+&+&+&+&-\\ &-&+&+&-&+&-\\ &+&+&+&-&+&-\\ &+&+&+&-&-&-&\\
&-&+&+&-&-&-\\ &-&+&+&+&-&-\\ &+&-&-&-&+&+\\ &+&-&-&+&+&+&\\
&-&-&-&+&+&+\\ &-&-&-&+&-&+\\ &+&-&-&+&-&+\\ &+&-&-&-&-&+&\\
&+&-&-&-&+&-\\ &+&+&-&-&+&-\\ &+&+&-&-&-&+\\ &+&+&-&+&-&+&\\
&-&+&-&+&-&+\\ &-&+&-&+&-&-\\ &+&+&-&+&-&-\\ &+&+&-&-&-&-&\}
\end{matrix}
\end{equation*}
}

The tope committee {\tiny
\begin{equation*}
\begin{matrix} \{&+&+&+&-&+&-&\\ &-&+&+&+&-&+&\\ &+&-&-&+&+&+&\\
&+&+&-&+&-&+&\\ &+&-&+&-&+&+&\\ &-&+&+&+&+&-&\}\\
\end{matrix}
\end{equation*}
} \hspace{-1.4mm}for $\mathcal{N}^2$, of even cardinality, is not
minimal; indeed, it splits up into a disjoint union of two
critical committees, {\tiny
\begin{equation*}
\begin{matrix} \{&+&+&+&-&+&-&\\ &-&+&+&+&-&+&\\ &+&-&-&+&+&+&\}
\end{matrix}\ \ \ \ \ \text{\normalsize and}\ \ \ \ \
\begin{matrix} \{&+&+&-&+&-&+&\\ &+&-&+&-&+&+&\\ &-&+&+&+&+&-&\}&\ ,
\end{matrix}
\end{equation*}
} \hspace{-1.5mm}of minimal cardinality, because
$\mathcal{N}^2$ is not acyclic.

As the topes of the committee {\tiny
\begin{equation*}
\begin{matrix} \{&+&+&+&-&+&-&&\\ &-&+&+&+&-&+&&\\ &+&+&-&+&-&+&&\\
&+&-&-&+&+&+&&\\ &+&-&+&-&+&+&\}&
\end{matrix}
\end{equation*}
} \hspace{-1.5mm}have the positive parts which are maximal
with respect to inclusion, it is, intuitively, of `higher quality' than the committee {\tiny
\begin{equation*}
\begin{matrix} \{&-&+&+&-&+&-&&\\ &-&+&+&+&-&-&&\\ &+&+&-&+&-&+&&\\
&+&-&-&+&+&+&&\\ &+&-&+&-&+&+&\}&\ .
\end{matrix}
\end{equation*}
}
\end{example}

\subsection{Committee Decision Rules}

Let $\mathcal{S}$ be a training set on a ground set~$E$. Denote by $\mathcal{M}$ the reorientation
\begin{equation*}
\mathcal{M}:={}_{-\lambda^{-1}(-)}\mathcal{S}\ .
\end{equation*}

Let $\widetilde{\mathcal{M}}:=(E\dot\cup\{g\},\widetilde{\mathcal{C}^{\ast}})$ denote a nontrivial single
element extension of $\mathcal{M}=(E,\mathcal{C}^{\ast})=(E,\mathcal{L})$ by a new pattern $g$, such that
$\widetilde{\mathcal{M}}$ is simple. Let $\sigma:\ \mathcal{C}^{\ast}\to\{-,0,+\}$ denote the corresponding
localization.

If $\mathcal{K}^{\ast}$ is a tope committee for $\mathcal{M}$,
then assign to every tope $K\in\mathcal{K}^{\ast}$ the set of
cocircuits
\begin{equation*}
\mathcal{C}_K{}^{\ast}:=\{D\in\mathcal{C}^{\ast}:\ \text{$D$
conforms to $K$}\}\ .
\end{equation*}

\begin{itemize}
\item
If the sets $\{\bigl(D,\sigma(D)\bigr):\
D\in\mathcal{C}_K{}^{\ast}\}$ are conformal, for all
$K\in\mathcal{K}^{\ast}$, then define a subset of topes
$\widetilde{\mathcal{K}^{\ast}}$ of $\widetilde{\mathcal{M}}$ in
the following way:
\begin{equation}
\label{eq:35} \widetilde{\mathcal{K}^{\ast}}:=
\Bigl\{\underset{D\in\mathcal{C}_{K}{}^{\ast}}\bigcirc
\bigl({}_{-\lambda^{-1}(-)}D,\sigma(D)\bigr):\
K\in\mathcal{K}^{\ast}\Bigr\}\ .
\end{equation}

The {\em committee decision rule\/} corresponding to
$\mathcal{K}^{\ast}$ is the mapping $\mathfrak{r}:\
E\dot\cup\{g\}\to\{-,0,+\}$ such that
\begin{equation}
\label{eq:18} \mathfrak{r}:\ f\mapsto\begin{cases}-\ ,&\text{if\ \
$\bigl|\{\widetilde{K}\in\widetilde{\mathcal{K}^{\ast}}:\
\widetilde{K}(f)=-\}\bigr|>
\bigl|\{\widetilde{K}\in\widetilde{\mathcal{K}^{\ast}}:\
\widetilde{K}(f)=+\}\bigr|$}\ \ ,\\ +\ ,&\text{if\ \
$\bigl|\{\widetilde{K}\in\widetilde{\mathcal{K}^{\ast}}:\
\widetilde{K}(f)=-\}\bigr|<
\bigl|\{\widetilde{K}\in\widetilde{\mathcal{K}^{\ast}}:\
\widetilde{K}(f)=+\}\bigr|$}\ \ ,\\ 0\ ,&\text{otherwise}\ .
\end{cases}
\end{equation}

\item
If $\{\bigl(D,\sigma(D)\bigr):\ D\in\mathcal{C}_K{}^{\ast}\}$ is
not conformal, for some $K\in\mathcal{K}^{\ast}$, then the {\em
committee decision rule\/} $\mathfrak{r}$ corresponding to
$\mathcal{K}^{\ast}$ is defined by
\begin{equation*}
\mathfrak{r}:\ e\mapsto\begin{cases}-\ ,&\text{if\ \
$\bigl|\{K\in\mathcal{K}^{\ast}:\
\left({}_{-\lambda^{-1}(-)}K\right)\!(e)=-\}\bigr|$}\\
\quad&\text{\quad\quad\quad $> \bigl|\{K\in\mathcal{K}^{\ast}:\
\left({}_{-\lambda^{-1}(-)}K\right)\!(e)=+\}\bigr|$}\ \ ,\\ +\
,&\text{if\ \ $\bigl|\{K\in\mathcal{K}^{\ast}:\
\left({}_{-\lambda^{-1}(-)}K\right)\!(e)=-\}\bigr|$}\\
\quad&\text{\quad\quad\quad $< \bigl|\{K\in\mathcal{K}^{\ast}:\
\left({}_{-\lambda^{-1}(-)}K\right)\!(e)=+\}\bigr|$}\ \ ,
\end{cases}
\end{equation*}
for all $e\in E$. By convention, $\mathfrak{r}:\ g\mapsto 0$.
\end{itemize}

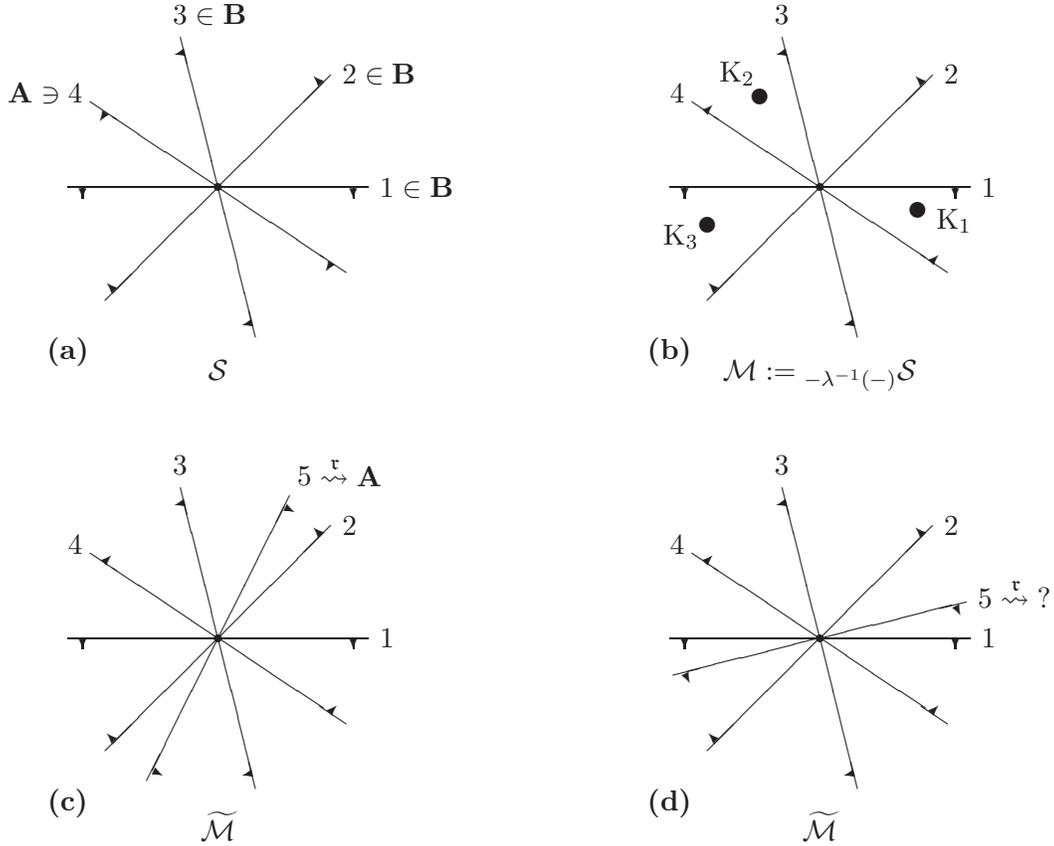
\begin{figure}[ht]
\begin{picture}(10,12)(3,0)

\put(4,9){\circle * {0.1}}

\put(4,9){\line(1,0){2}} \put(4,9){\line(-1,0){2}}
\put(4,9){\line(1,1){1.5}} \put(4,9){\line(-1,-1){1.5}}
\put(4,9){\line(1,-4){0.5}} \put(4,9){\line(-1,4){0.5}}
\put(4,9){\line(3,-2){1.7}} \put(4,9){\line(-3,2){1.7}}

\put(2.2,9){\vector(0,-1){0.15}} \put(5.8,9){\vector(0,-1){0.15}}
\put(2.65,7.6){\vector(-1,1){0.13}}
\put(5.4,10.35){\vector(-1,1){0.13}}
\put(4.45,7.2){\vector(-4,-1){0.13}}
\put(3.55,10.8){\vector(-4,-1){0.13}}
\put(5.58,8.1){\vector(-2,-3){0.15}}
\put(2.58,10.1){\vector(-2,-3){0.15}}

\put(6.15,9){\makebox(0,0)[l]{$1\in\mathbf{B}$}}
\put(5.65,10.5){\makebox(0,0)[l]{$2\in\mathbf{B}$}}
\put(3.4,11.3){\makebox(0,0)[l]{$3\in\mathbf{B}$}}
\put(2.2,10.25){\makebox(0,0)[r]{$\mathbf{A}\ni 4$}}

\put(4,6.7){\makebox(0,0)[t]{$\mathcal{S}$}}
\put(2,7){\makebox(0,0)[t]{\bf{(a)}}}


\put(12,9){\circle * {0.1}}

\put(12,9){\line(1,0){2}} \put(12,9){\line(-1,0){2}}
\put(12,9){\line(1,1){1.5}} \put(12,9){\line(-1,-1){1.5}}
\put(12,9){\line(1,-4){0.5}} \put(12,9){\line(-1,4){0.5}}
\put(12,9){\line(3,-2){1.7}} \put(12,9){\line(-3,2){1.7}}

\put(10.2,9){\vector(0,-1){0.15}}
\put(13.8,9){\vector(0,-1){0.15}}
\put(10.65,7.6){\vector(-1,1){0.13}}
\put(13.4,10.35){\vector(-1,1){0.13}}
\put(12.45,7.2){\vector(-4,-1){0.13}}
\put(11.55,10.8){\vector(-4,-1){0.13}}
\put(13.43,7.95){\vector(1,1){0.15}}
\put(10.43,9.95){\vector(1,1){0.15}}

\put(14.15,9){\makebox(0,0)[l]{$1$}}
\put(13.65,10.5){\makebox(0,0)[l]{$2$}}
\put(11.4,11.3){\makebox(0,0)[l]{$3$}}
\put(10.2,10.25){\makebox(0,0)[r]{$4$}}

\put(12,6.7){\makebox(0,0)[t]{$\mathcal{M}:={}_{-\lambda^{-1}(-)}\mathcal{S}$}}

\put(10.5,8.5){\circle * {0.2}}
\put(10.15,8.5){\makebox(0,0)[t]{$\mathrm{K}_3$}}

\put(13.3,8.7){\circle * {0.2}}
\put(13.55,8.55){\makebox(0,0)[l]{$\mathrm{K}_1$}}

\put(11.2,10.2){\circle * {0.2}}
\put(10.9,10.35){\makebox(0,0)[b]{$\mathrm{K}_2$}}

\put(10,7){\makebox(0,0)[t]{\bf{(b)}}}


\put(4,3){\circle * {0.1}}

\put(4,3){\line(1,0){2}} \put(4,3){\line(-1,0){2}}
\put(4,3){\line(1,1){1.5}} \put(4,3){\line(-1,-1){1.5}}
\put(4,3){\line(1,-4){0.5}} \put(4,3){\line(-1,4){0.5}}
\put(4,3){\line(3,-2){1.7}} \put(4,3){\line(-3,2){1.7}}

\put(4,3){\line(1,2){0.95}} \put(4,3){\line(-1,-2){0.95}}

\put(2.2,3){\vector(0,-1){0.15}} \put(5.8,3){\vector(0,-1){0.15}}
\put(2.65,1.6){\vector(-1,1){0.13}}
\put(5.4,4.35){\vector(-1,1){0.13}}
\put(4.45,1.2){\vector(-4,-1){0.13}}
\put(3.55,4.8){\vector(-4,-1){0.13}}
\put(5.43,1.95){\vector(1,1){0.15}}
\put(2.43,3.95){\vector(1,1){0.15}}

\put(4.87,4.75){\vector(2,-1){0.15}}
\put(3.12,1.25){\vector(2,-1){0.15}}

\put(6.15,3){\makebox(0,0)[l]{$1$}}
\put(5.65,4.5){\makebox(0,0)[l]{$2$}}
\put(3.4,5.3){\makebox(0,0)[l]{$3$}}
\put(2.2,4.25){\makebox(0,0)[r]{$4$}}

\put(5.15,5.15){\makebox(0,0)[c]{$5$}}
\put(5.35,5.22){\makebox(0,0)[l]{$\overset{\mathfrak{r}}\rightsquigarrow\mathbf{A}$}}

\put(4,0.7){\makebox(0,0)[t]{$\widetilde{\mathcal{M}}$}}
\put(2,1){\makebox(0,0)[t]{\bf{(c)}}}


\put(12,3){\circle * {0.1}}

\put(12,3){\line(1,0){2}} \put(12,3){\line(-1,0){2}}
\put(12,3){\line(1,1){1.5}} \put(12,3){\line(-1,-1){1.5}}
\put(12,3){\line(1,-4){0.5}} \put(12,3){\line(-1,4){0.5}}
\put(12,3){\line(3,-2){1.7}} \put(12,3){\line(-3,2){1.7}}

\put(12,3){\line(4,1){1.95}} \put(12,3){\line(-4,-1){1.95}}

\put(10.2,3){\vector(0,-1){0.15}}
\put(13.8,3){\vector(0,-1){0.15}}
\put(10.65,1.6){\vector(-1,1){0.13}}
\put(13.4,4.35){\vector(-1,1){0.13}}
\put(12.45,1.2){\vector(-4,-1){0.13}}
\put(11.55,4.8){\vector(-4,-1){0.13}}
\put(13.43,1.95){\vector(1,1){0.15}}
\put(10.43,3.95){\vector(1,1){0.15}}

\put(13.72,3.6){\vector(1,-2){0.15}}
\put(10.12,2.7){\vector(1,-2){0.15}}

\put(14.15,3){\makebox(0,0)[l]{$1$}}
\put(13.65,4.5){\makebox(0,0)[l]{$2$}}
\put(11.4,5.3){\makebox(0,0)[l]{$3$}}
\put(10.2,4.25){\makebox(0,0)[r]{$4$}}

\put(14.1,3.6){\makebox(0,0)[l]{$5\overset{\mathfrak{r}}\rightsquigarrow\bf{?}$}}

\put(12,0.7){\makebox(0,0)[t]{$\widetilde{\mathcal{M}}$}}
\put(10,1){\makebox(0,0)[t]{\bf{(d)}}}
\end{picture}
\caption{(a): A central line arrangement that realizes a rank $2$
training set $\mathcal{S}$. The positive halfplanes of
$\mathbb{R}^2$ are marked by arrows; (b): A realization of the
reorientation $\mathcal{M}:={}_{-\lambda^{-1}(-)}\mathcal{S}$. The
set $\{\mathrm{K}_1,\mathrm{K}_2,\mathrm{K}_3\}$ of regions marked
by discs corresponds to a tope committee for $\mathcal{M}$; (c):
The new pattern $5$ is classified as an element of the class
$\mathbf{A}$; (d): The new pattern $5$ is not classified.}
\label{fig:9}
\end{figure}

\begin{example}
Figure~{\rm\ref{fig:9}(a)} depicts a realization of a rank $2$
training set\/ $\mathcal{S}$ on the ground set $E_4$. A
realization of its reorientation
$\mathcal{M}:={}_{-\lambda^{-1}(-)}\mathcal{S}={}_{-4}\mathcal{S}$
is shown in Figure~{\rm\ref{fig:9}(b)}. The set of topes {\small
\begin{equation*}
\mathcal{K}^{\ast}:=
\begin{matrix}\{&K_1:=&+&-&-&+&\\ &K_2:=&-&+&+&+&\\ &K_3:=&+&+&+&-&\}
\end{matrix}
\end{equation*}
} \hspace{-1.2mm}is a committee for $\mathcal{M}$.

\begin{itemize}
\item Let $\widetilde{\mathcal{M}}$ be a nontrivial single element
extension of $\mathcal{M}$ by the pattern~$5$, as shown in
Figure~{\rm\ref{fig:9}(c)}.

Each of the sets of cocircuits {\small
\begin{align*}
\{\bigl(D,\sigma(D)\bigr):\
D\in\mathcal{C}_{K_1}{}^{\ast}\}&=\begin{matrix}\{&0&-&-&+&+&\\
&+&-&-&0&+&\}\end{matrix}\ ,\\ \{\bigl(D,\sigma(D)\bigr):\
D\in\mathcal{C}_{K_2}{}^{\ast}\}&=\begin{matrix}\{&-&+&0&+&-&\\
&-&+&+&0&-&\}\end{matrix}\ ,\\ \{\bigl(D,\sigma(D)\bigr):\
D\in\mathcal{C}_{K_3}{}^{\ast}\}&=\begin{matrix}\{&0&+&+&-&-&\\
&+&0&+&-&-&\}\end{matrix}
\end{align*}
} \hspace{-1.5mm}is conformal. The set of topes
$\widetilde{\mathcal{K}^{\ast}}$, defined by~{\rm(\ref{eq:35})},
is {\small
\begin{equation*}
\widetilde{\mathcal{K}^{\ast}}=
\begin{matrix}\{&+&-&-&-&+&&\\ &-&+&+&-&-&&\\ &+&+&+&+&-&\}&.
\end{matrix}
\end{equation*}
} \hspace{-1.5mm}Therefore the decision rule $\mathfrak{r}$, corresponding to $\mathcal{K}^{\ast}$ and
defined by~{\rm(\ref{eq:18})}, recognizes the pattern~$5$ as an element of the class $\mathbf{A}$.
\item
If $\widetilde{\mathcal{M}}$ is a nontrivial single element
extension of $\mathcal{M}$ by the pattern $5$, as shown in
Figure~{\rm\ref{fig:9}(d)}, then the set of cocircuits {\small
\begin{equation*}
\{\bigl(D,\sigma(D)\bigr):\
D\in\mathcal{C}_{K_3}{}^{\ast}\}=\begin{matrix}\{&0&+&+&-&-&\\
&+&0&+&-&+&\}\end{matrix}
\end{equation*}
} \hspace{-1.4mm}is not conformal. As a consequence,
$\mathfrak{r}(g):=0$, that is, the decision rule $\mathfrak{r}$
corresponding to $\mathcal{K}^{\ast}$ leaves the pattern $5$
unclassified.
\end{itemize}
\end{example}

\section{The Existence of a Tope Committee: Reorientations}
\label{section:3} In this section we show that every simple oriented matroid has a tope committee. In fact, a
critical tope committee for such an oriented matroid can be built based on information on an arbitrary
maximal chain in the tope poset.

Our argument relies on the mechanism of consecutive reorientations
of an initial acyclic oriented matroid.

The main construction which we make use of is a direct
generalization of centrally-symmetric cycles of regions in a
central hyperplane arrangement from~\cite{AK-IEEE}:

\begin{remark}
\label{p:8} Let $\mathcal{M}:=(E_m,\mathcal{L})=(E_m,\mathcal{T})$
be a simple oriented matroid. Let
$\pmb{R}:=(T^0,T^1,\ldots,T^{2m-1},T^0)$ be a cycle in the tope
graph $\mathcal{T}(\mathcal{L})$ such that
\begin{equation*}
T^{k+m}=-T^k\ ,\ \ \ 0\leq k\leq m-1\ .
\end{equation*}
We call such a cycle\/ {\em symmetric}.
\begin{itemize}
\item[\rm(i)]
For every $e\in E_m$ the set of topes
\begin{equation*}
\{T\in\mathfrak{V}(\pmb{R}):\ T(e)=+\}
\end{equation*}
is the set of vertices of a path of length $m-1$ in $\pmb{R}$; if a tope $T$ is an endvertex of the path then
the other endvertex is the tope ${}_{-e}(-T)$.
\item[\rm(ii)] Let $(T^{k_1},T^{k_2},T^{k_3})$ be a $2$-path in
$\pmb{R}$. We have
\begin{equation}
\label{eq:19} T^{k_2} \in\bmax^+\bigl(\mathfrak{V}(\pmb{R})\bigr)
\end{equation}
iff $(T^{k_1})^+\subsetneqq(T^{k_2})^+\supsetneqq(T^{k_3})^+$ or, equivalently, $|(T^{k_1})^+|=|(T^{k_3})^+|$
$=|(T^{k_2})^+|-1$.

In other words, let $\{f\}:=\mathbf{S}(T^{k_1},T^{k_2})$ and $\{g\}:=\mathbf{S}(T^{k_2},T^{k_3})$; then
inclusion~{\rm(\ref{eq:19})} holds iff\/ $T^{k_2}(f)=T^{k_2}(g)=+$.
\end{itemize}
\end{remark}

\subsection{Rank $2$}

In the theory of oriented matroids the rank $2$ case is
instructive.

\begin{lemma}
\label{p:6} Let $\mathcal{N}^0:=(E_m,\mathcal{T}^0)$ be a simple
acyclic oriented matroid of rank $2$, on the ground set $E_m$,
with set of topes $\mathcal{T}^0$.

Let $(j_1,\ldots,j_s)$ be a nonempty sequence of integers with $j_i\in E_m$, $1\leq i$ $\leq s$. Define the
reorientation $\mathcal{N}^{i}:=(E_m,\mathcal{T}^i):={}_{-j_i}\mathcal{N}^{i-1}$ whose set of topes is
denoted by $\mathcal{T}^i$.

The reorientation $\mathcal{N}^{s}$ of $\mathcal{N}^0$ has a
critical tope committee.
\end{lemma}

Three types of transformations of committees are carried out by Algorithm~\ref{alg:1} that underlies the
proof below. We illustrate the transformations by considering several rank $2$ oriented matroids represented
by central line arrangements in the plane. The regions corresponding to the topes from committees will be
marked in figures by discs or circles. Let $\mathcal{K}_0{}^{\ast}:=\{\mathrm{T}^{(+)}\}$, and let
$\mathcal{K}_i{}^{\ast}$ denote the tope committee built by Algorithm~\ref{alg:1} for the reorientation
$\mathcal{N}^i$.

\begin{itemize}
\item[1)] If there is a tope $K\in\mathcal{K}_{i-1}{}^{\ast}$
such that
\begin{equation*}
K(j_i)=+\ \ \ \text{and}\ \ \  {}_{-j_i}(-K)\in\mathcal{T}^{i-1}
\end{equation*}
(or, equivalently, $K(j_i)=+$, and there is a subtope $H <
\!\!\!\cdot\ K$ such that $\mathbf{z}(H)=\{j_i\}$), but there is
no tope $S$ in $\mathcal{K}_{i-1}{}^{\ast}$ such that $S(j_i)=+$
and ${}_{-j_i}S=-K$, then the set
\begin{equation*}
\mathcal{K}_i{}^{\ast}:={}_{-j_i}(\mathcal{K}_{i-1}{}^{\ast}-\{K\})\
\dot\cup\ \{K\}
\end{equation*}
is a tope committee for $\mathcal{N}^i$; see Figure~\ref{fig:3}.

\begin{figure}[ht]
\begin{picture}(10,5)(3,1.5)

\put(4,4){\circle * {0.1}}

\put(4,4){\line(1,0){2}} \put(4,4){\line(-1,0){2}}
\put(4,4){\line(1,1){1.5}} \put(4,4){\line(-1,-1){1.5}}
\put(4,4){\line(1,-4){0.5}} \put(4,4){\line(-1,4){0.5}}
\put(4,4){\line(2,-3){1.15}} \put(4,4){\line(-2,3){1.15}}
\put(4,4){\line(3,-2){1.7}} \put(4,4){\line(-3,2){1.7}}

\put(2.2,4){\vector(0,-1){0.15}} \put(5.8,4){\vector(0,-1){0.15}}
\put(2.55,2.6){\vector(1,-1){0.13}}
\put(5.35,5.4){\vector(1,-1){0.13}}
\put(4.45,2.2){\vector(-4,-1){0.13}}
\put(3.55,5.8){\vector(-4,-1){0.13}}
\put(4.97,2.5){\vector(3,2){0.15}}
\put(2.95,5.5){\vector(3,2){0.15}}
\put(5.58,3.1){\vector(-2,-3){0.15}}
\put(2.58,5.1){\vector(-2,-3){0.15}}

\put(3.75,3){\circle * {0.2}} \put(4.8,3.2){\circle * {0.2}}
\put(3.55,5){\circle * {0.2}}

\put(5.8,2.8){\makebox(0,0)[l]{$j_i$}}

\put(5.3,2.5){\makebox(0,0)[l]{$K$}}

\put(4,1.7){\makebox(0,0)[t]{$\mathcal{N}^{i-1}$}}


\put(12,4){\circle * {0.1}}

\put(12,4){\line(1,0){2}} \put(12,4){\line(-1,0){2}}
\put(12,4){\line(1,1){1.5}} \put(12,4){\line(-1,-1){1.5}}
\put(12,4){\line(1,-4){0.5}} \put(12,4){\line(-1,4){0.5}}
\put(12,4){\line(2,-3){1.15}} \put(12,4){\line(-2,3){1.15}}
\put(12,4){\line(3,-2){1.7}} \put(12,4){\line(-3,2){1.7}}

\put(10.2,4){\vector(0,-1){0.15}}
\put(13.8,4){\vector(0,-1){0.15}}
\put(10.55,2.6){\vector(1,-1){0.13}}
\put(13.35,5.4){\vector(1,-1){0.13}}
\put(12.45,2.2){\vector(-4,-1){0.13}}
\put(11.55,5.8){\vector(-4,-1){0.13}}
\put(12.97,2.5){\vector(3,2){0.15}}
\put(10.95,5.5){\vector(3,2){0.15}}
\put(13.4,2.9){\vector(2,3){0.15}}
\put(10.42,4.9){\vector(2,3){0.15}}

\put(11.75,3){\circle * {0.2}} \put(13.35,3.65){\circle  {0.18}}
\put(11.55,5){\circle * {0.2}}

\put(13.8,2.8){\makebox(0,0)[l]{$j_i$}}
\put(12,1.7){\makebox(0,0)[t]{$\mathcal{N}^i$}}

\end{picture}
\caption{A transformation of a tope committee under a reorientation:
$\mathcal{K}_i{}^{\ast}:={}_{-j_i}(\mathcal{K}_{i-1}{}^{\ast}-\{K\})\ \dot\cup\ \{K\}$; here
$|\mathcal{K}_i{}^{\ast}|=|\mathcal{K}_{i-1}{}^{\ast}|$.} \label{fig:3}
\end{figure}
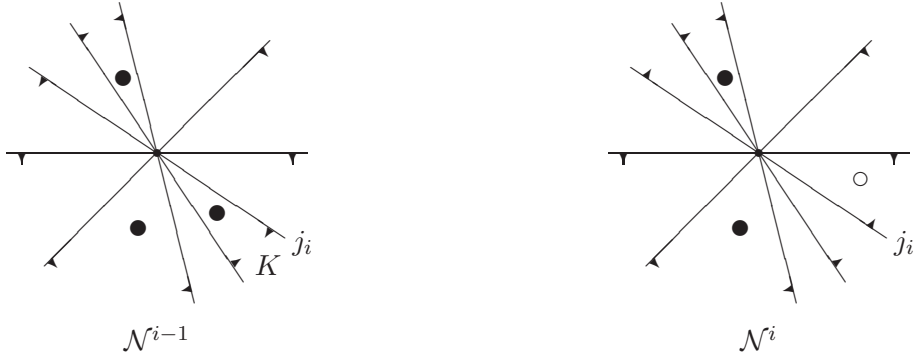

\item[2)] If there is no tope $K$ in $\mathcal{K}_{i-1}{}^{\ast}$
such that $K(j_i)=+$ and ${}_{-j_i}(-K)$ $\in\mathcal{T}^{i-1}$, then pick the pair of topes
$\{T'',T'''\}\subset\mathcal{T}^{i-1}$ such that $T''(j_i)=T'''(j_i)=-$ and ${}_{-j_i}T''=-T'''$. The set
\begin{equation*}
\mathcal{K}_i{}^{\ast}:={}_{-j_i}\bigl(\mathcal{K}_{i-1}{}^{\ast}\
\dot{\cup}\ \{T'',\ T'''\}\bigr)
\end{equation*} is a tope committee for
$\mathcal{N}^i$; see Figure~\ref{fig:1}.

\begin{figure}[ht]
\begin{picture}(10,5)(3,1.5)

\put(4,4){\circle * {0.1}}

\put(4,4){\line(1,0){2}} \put(4,4){\line(-1,0){2}}
\put(4,4){\line(1,1){1.5}} \put(4,4){\line(-1,-1){1.5}}
\put(4,4){\line(1,-4){0.5}} \put(4,4){\line(-1,4){0.5}}
\put(4,4){\line(2,-3){1.15}} \put(4,4){\line(-2,3){1.15}}
\put(4,4){\line(3,-2){1.7}} \put(4,4){\line(-3,2){1.7}}

\put(2.2,4){\vector(0,-1){0.15}} \put(5.8,4){\vector(0,-1){0.15}}
\put(2.55,2.6){\vector(1,-1){0.13}}
\put(5.35,5.4){\vector(1,-1){0.13}}
\put(4.45,2.2){\vector(-4,-1){0.13}}
\put(3.55,5.8){\vector(-4,-1){0.13}}
\put(4.97,2.5){\vector(3,2){0.15}}
\put(2.95,5.5){\vector(3,2){0.15}}
\put(5.58,3.1){\vector(-2,-3){0.15}}
\put(2.58,5.1){\vector(-2,-3){0.15}}

\put(3.75,3){\circle * {0.2}} \put(4.8,3.2){\circle * {0.2}}
\put(3.55,5){\circle * {0.2}}

\put(6.1,4){\makebox(0,0)[l]{$j_i$}}

\put(2.5,4.5){\makebox(0,0)[r]{$T''$}}
\put(5.4,4.5){\makebox(0,0)[l]{$T'''$}}

\put(4,1.7){\makebox(0,0)[t]{$\mathcal{N}^{i-1}$}}


\put(12,4){\circle * {0.1}}

\put(12,4){\line(1,0){2}} \put(12,4){\line(-1,0){2}}
\put(12,4){\line(1,1){1.5}} \put(12,4){\line(-1,-1){1.5}}
\put(12,4){\line(1,-4){0.5}} \put(12,4){\line(-1,4){0.5}}
\put(12,4){\line(2,-3){1.15}} \put(12,4){\line(-2,3){1.15}}
\put(12,4){\line(3,-2){1.7}} \put(12,4){\line(-3,2){1.7}}

\put(10.2,4){\vector(0,1){0.15}} \put(13.8,4){\vector(0,1){0.15}}
\put(10.55,2.6){\vector(1,-1){0.13}}
\put(13.35,5.4){\vector(1,-1){0.13}}
\put(12.45,2.2){\vector(-4,-1){0.13}}
\put(11.55,5.8){\vector(-4,-1){0.13}}
\put(12.97,2.5){\vector(3,2){0.15}}
\put(10.95,5.5){\vector(3,2){0.15}}
\put(13.58,3.1){\vector(-2,-3){0.15}}
\put(10.58,5.1){\vector(-2,-3){0.15}}

\put(11.75,3){\circle * {0.2}} \put(12.8,3.2){\circle * {0.2}}
\put(11.55,5){\circle * {0.2}}

\put(10.8,4.25){\circle {0.18}} \put(13,4.5){\circle {0.18}}

\put(14.1,4){\makebox(0,0)[l]{$j_i$}}
\put(12,1.7){\makebox(0,0)[t]{$\mathcal{N}^i$}}

\end{picture}
\caption{A transformation of a tope committee under a reorientation:
$\mathcal{K}_i{}^{\ast}:={}_{-j_i}\bigl(\mathcal{K}_{i-1}{}^{\ast}\ \dot{\cup}\ \{T'',\ T'''\}\bigr)$; here
$|\mathcal{K}_i{}^{\ast}|=|\mathcal{K}_{i-1}{}^{\ast}|+2$.} \label{fig:1}
\end{figure}

\item[3)] If there is a tope $K\in\mathcal{K}_{i-1}{}^{\ast}$
such that $K(j_i)=+$ and ${}_{-j_i}(-K)$ $\in\mathcal{T}^{i-1}$, and if there is a tope
$S\in\mathcal{K}_{i-1}{}^{\ast}$ such that ${}_{-j_i}S=-K$, then the set
\begin{equation*}
\mathcal{K}_i{}^{\ast}:={}_{-j_i}\bigl(\mathcal{K}_{i-1}{}^{\ast}-\{K,S\}\bigr)
\end{equation*}
is a tope committee for $\mathcal{N}^i$; see Figure~\ref{fig:2}.

\begin{figure}[ht]
\begin{picture}(10,5)(3,1.5)

\put(4,4){\circle * {0.1}}

\put(4,4){\line(1,0){2}} \put(4,4){\line(-1,0){2}}
\put(4,4){\line(1,1){1.5}} \put(4,4){\line(-1,-1){1.5}}
\put(4,4){\line(1,-4){0.5}} \put(4,4){\line(-1,4){0.5}}
\put(4,4){\line(2,-3){1.15}} \put(4,4){\line(-2,3){1.15}}
\put(4,4){\line(3,-2){1.7}} \put(4,4){\line(-3,2){1.7}}

\put(2.2,4){\vector(0,1){0.15}} \put(5.8,4){\vector(0,1){0.15}}
\put(2.55,2.6){\vector(1,-1){0.13}}
\put(5.35,5.4){\vector(1,-1){0.13}}
\put(4.45,2.2){\vector(-4,-1){0.13}}
\put(3.55,5.8){\vector(-4,-1){0.13}}
\put(4.97,2.5){\vector(3,2){0.15}}
\put(2.95,5.5){\vector(3,2){0.15}}
\put(5.58,3.1){\vector(-2,-3){0.15}}
\put(2.58,5.1){\vector(-2,-3){0.15}}

\put(3.75,3){\circle * {0.2}} \put(4.8,3.2){\circle * {0.2}}
\put(3.55,5){\circle * {0.2}}

\put(2.8,4.25){\circle * {0.2}} \put(5,4.5){\circle * {0.2}}

\put(2.5,2.3){\makebox(0,0)[r]{$j_i$}}

\put(3.5,2.3){\makebox(0,0)[r]{$S$}}
\put(6,4.8){\makebox(0,0)[r]{$K$}}

\put(4,1.7){\makebox(0,0)[t]{$\mathcal{N}^{i-1}$}}


\put(12,4){\circle * {0.1}}

\put(12,4){\line(1,0){2}} \put(12,4){\line(-1,0){2}}
\put(12,4){\line(1,1){1.5}} \put(12,4){\line(-1,-1){1.5}}
\put(12,4){\line(1,-4){0.5}} \put(12,4){\line(-1,4){0.5}}
\put(12,4){\line(2,-3){1.15}} \put(12,4){\line(-2,3){1.15}}
\put(12,4){\line(3,-2){1.7}} \put(12,4){\line(-3,2){1.7}}

\put(10.2,4){\vector(0,1){0.15}} \put(13.8,4){\vector(0,1){0.15}}
\put(10.65,2.6){\vector(-1,1){0.13}}
\put(13.45,5.4){\vector(-1,1){0.13}}
\put(12.45,2.2){\vector(-4,-1){0.13}}
\put(11.55,5.8){\vector(-4,-1){0.13}}
\put(12.97,2.5){\vector(3,2){0.15}}
\put(10.95,5.5){\vector(3,2){0.15}}
\put(13.58,3.1){\vector(-2,-3){0.15}}
\put(10.58,5.1){\vector(-2,-3){0.15}}

\put(12.8,3.2){\circle * {0.2}} \put(11.55,5){\circle * {0.2}}
\put(10.8,4.25){\circle * {0.2}}

\put(10.5,2.3){\makebox(0,0)[r]{$j_i$}}
\put(12,1.7){\makebox(0,0)[t]{$\mathcal{N}^i$}}

\end{picture}
\caption{A transformation of a tope committee under a reorientation:
$\mathcal{K}_i{}^{\ast}:={}_{-j_i}\bigl(\mathcal{K}_{i-1}{}^{\ast}-\{K,S\}\bigr)$; here
$|\mathcal{K}_i{}^{\ast}|=|\mathcal{K}_{i-1}{}^{\ast}|-2$.} \label{fig:2}
\end{figure}

\end{itemize}

\begin{proof}
The set $\mathcal{K}_s{}^{\ast}\subset\mathcal{T}^s$ built by
means of Algorithm~\ref{alg:1} is a tope committee for
$\mathcal{N}^s$:

{\sl Claim {\rm 1}.} {\em For any $i$, $1\leq i\leq s$,
$\mathcal{K}_i{}^{\ast}$ is of odd cardinality, and for any $e\in
E_m$, it holds
\begin{equation}
\label{eq:3} |\{K\in\mathcal{K}_i{}^{\ast}:\ K(e)=+\}|
=\left\lceil\tfrac{|\mathcal{K}_i{}^{\ast}|}{2}\right\rceil\ .
\end{equation}
}

{\sl Claim {\rm 2}.} {\em For any $i$, $1\leq i\leq s$, we have
\begin{equation}
\label{eq:4} \mathcal{K}_i{}^{\ast}=\bigl\{K\in\mathcal{T}^i:\
T\in\mathcal{T}^i,\ \mathbf{S}(K,T)=\{e\}\ \Longrightarrow\ K(e)=+
\bigr\}\ .
\end{equation}
As a consequence, the committee $\mathcal{K}_i{}^{\ast}$ is
minimal.}

$\vartriangleright$ Indeed, let $i=1$.
\begin{itemize}
\item[$\bullet$]
If conditions (\ref{eq:1}) hold for $K:=\mathrm{T}^{(+)}$ in
$\mathcal{N}^{i-1}$, then pick the tope $T\in\mathcal{T}^{i-1}$
such that $\mathbf{S}(K,T)=\{j_i\}$ in $\mathcal{N}^{i-1}$; the
one-element set
$\mathcal{K}_i{}^{\ast}:=\left\{{}_{-j_i}T\right\}=\{K\}=\{\mathrm{T}^{(+)}\}$,
formed at Step~{\it 10}\/ of the algorithm, is a tope committee
for $\mathcal{N}^i$.

\item[$\bullet$]
If (\ref{eq:1}) do not hold for $K:=\mathrm{T}^{(+)}$, then the
three-element set of covectors
$\mathcal{K}_i{}^{\ast}:=\{K',K'',K'''\}$, built at Step~{\it
14}\/ of the algorithm, where
\begin{equation*}
K':={}_{-j_i}K\ ,\ \ \ \ K'':={}_{-j_i}T''\ ,\ \ \ \
K''':={}_{-j_i}T'''\ ,
\end{equation*}
is a tope committee for $\mathcal{N}^i$ because we have {\small
\begin{equation*}
\begin{matrix}\mathcal{K}_i{}^{\ast}=&\Bigl\{&K'=&+&
\ldots&\underset{\substack{\uparrow\\(j_i-1)}}{+}&
\underset{\substack{\uparrow\\j_i}}{-}&
\underset{\substack{\uparrow\\(j_i+1)}}{+}&\ldots&+&\phantom{\Bigr\}}
\\ \phantom{\mathcal{K}^{\ast}:=}&\phantom{\bigl\{}&
K''=&?&\ldots&?&+&?&\ldots& ?&\phantom{\Bigr\}}\\
\phantom{\mathcal{K}^{\ast}:=}&\phantom{\bigl\{}&
K'''=&-\bigl(K''(1)\bigr)&\ldots&-\bigl(K''(j-1)\bigr)&+&-\bigl(K''(j+1)\bigr)&
\ldots&-\bigl(K''(m)\bigr)&\Bigr\}\ ,\end{matrix}
\end{equation*}
} \hspace{-1.2mm}that is, for every $e\in E_m$, it holds
$|\{K\in\mathcal{K}_i{}^{\ast}:\
K(e)=+\}|=2=\left\lceil\tfrac{|\mathcal{K}_i{}^{\ast}|}{2}\right\rceil$.
\end{itemize}

Note that~(\ref{eq:3}) and~(\ref{eq:4}) hold for $i=1$.

$\quad$

Now, let $i>1$.

\begin{itemize}
\item[$\bullet$]
If there is no tope $K$ in $\mathcal{K}_{i-1}{}^{\ast}$ such
that~(\ref{eq:1}) hold, then for the topes ${}_{-j_i}T''$ and
${}_{-j_i}T'''$ added to the set $\mathcal{K}_i{}^{\ast}$ at
Step~{\it 14}\/ of the algorithm, we have
$\left({}_{-j_i}T''\right)(j_i)=\left({}_{-j_i}T'''\right)(j_i)=+$.
Assume that for the tope $T\in\mathcal{T}^{i-1}$ such that
$\mathbf{S}(T'',T)=\{f\}$ and $f\neq j_i$, it holds $T''(f)=-$.
Then the tope $K:=-T''$, with $K(j_i)=K(f)=+$, must belong to the
committee $\mathcal{K}_{i-1}{}^{\ast}$ and satisfy
conditions~(\ref{eq:1}), but this contradicts the negative
decision made at Step~{\it 06}. Hence $T''(e)=+$. For the tope
$Q\in\mathcal{T}^{i-1}$ such that $\mathbf{S}(T''',Q)=\{g\}$ and
$g\neq j_i$, we also have $T'''(g)=+$. As a result, (\ref{eq:4})
and (\ref{eq:3}) hold; hence $\mathcal{K}_i{}^{\ast}$ is a tope
committee for $\mathcal{N}^i$, of cardinality
$|\mathcal{K}_i{}^{\ast}|=|\mathcal{K}_{i-1}{}^{\ast}|+2$.

\item[$\bullet$]
If there are topes $K,S\in\mathcal{K}_{i-1}{}^{\ast}$ such that conditions~(\ref{eq:1}) and~(\ref{eq:2})
hold, then the algorithm excludes them from consideration at Steps~{\it 09\/} and~{\it 12}, and we obtain
$|\mathcal{K}_i{}^{\ast}|=|\mathcal{K}_{i-1}{}^{\ast}|-2$.

If there is a tope $K\in\mathcal{K}_{i-1}{}^{\ast}$ such
that~(\ref{eq:1}) hold, but there is no tope $S$ in
$\mathcal{K}_{i-1}{}^{\ast}$ satisfying~(\ref{eq:2}), then we have
$|\mathcal{K}_i{}^{\ast}|=|\mathcal{K}_{i-1}{}^{\ast}|$.

In any case, (\ref{eq:3}) holds, for all $e\in E_m$; therefore
$\mathcal{K}_i{}^{\ast}$ is a tope committee for $\mathcal{N}^i$.
\end{itemize}

If $\mathrm{T}^{(+)}\not\in\mathcal{K}_s{}^{\ast}$ then assume
that $\mathcal{K}_s{}^{\ast}$ is not minimal, that is, there is a
proper subset $\mathcal{Q}^{\ast}$ of the set
$\mathcal{K}_s{}^{\ast}$ such that
$\mathcal{K}_s{}^{\ast}-\mathcal{Q}^{\ast}$ is a committee for
$\mathcal{M}$. Since
$\mathrm{T}^{(+)}\not\in\mathcal{K}_s{}^{\ast}$, we have
$|\mathcal{K}_s{}^{\ast}-\mathcal{Q}^{\ast}|>1$.

Denote by $\pmb{R}:=(T^0,T^1,\ldots,T^{2m-1},T^0)$ the cycle which is the tope graph for $\mathcal{N}^s$, and
without loss of generality suppose that $T^0\in\mathcal{K}_s{}^{\ast}$ and $T^0$
\mbox{$\not\in\mathcal{Q}^{\ast}$}. Recall that $\mathcal{K}_s{}^{\ast}$ is precisely the set
$\bmax^+\bigl(\mathfrak{V}(\pmb{R})\bigr)=\bmax^+(\mathcal{T}^s)$, see Remark~\ref{p:8}(ii).

Let $\{g\}:=\mathbf{S}(T^0,T^1)$ and
$\{f\}:=\mathbf{S}(T^{2m-1},T^0)$; note that $f\neq g$. We have
$T^0(f)=T^0(g)=+$. For every $k$, $0<k<m$, we have $T^k(f)=+$; for
every $l$, $m<l<2m$, we have $T^l(g)=+$, see Remark~\ref{p:8}(i).

Claim~1 implies that for every $e\in\{f,g\}$ it holds $|\{Q\in\mathcal{Q}^{\ast}:\ Q(e)=$
$-\}|=|\{R\in\mathcal{Q}^{\ast}:\ R(e)=+\}|$, and we have
\begin{equation*}
|\{T^1,T^2,\ldots,T^{m-1}\}\cap\mathcal{Q}^{\ast}|=
|\{T^{m+1},T^{m+2},\ldots,T^{2m-1}\}\cap\mathcal{Q}^{\ast}|=
\frac{|\mathcal{Q}^{\ast}|}{2}\ ;
\end{equation*}
in particular, $\mathcal{Q}^{\ast}$ is of even cardinality.

Let $T^{k_2}$ be a tope from the set $\mathcal{Q}^{\ast}$, and
$(T^{k_1},T^{k_2},T^{k_3})$ a $2$-path in $\pmb{R}$. If
$\{p\}:=\mathbf{S}(T^{k_1},T^{k_2})$ and
$\{q\}:=\mathbf{S}(T^{k_2},T^{k_3})$, then there is an element
$h\in\{p,q\}$ such that
$|\{K\in\mathcal{K}_s{}^{\ast}-\mathcal{Q}^{\ast}:\
K(h)=+\}|=\left\lfloor\tfrac{|\mathcal{K}_s{}^{\ast}-\mathcal{Q}^{\ast}|}{2}\right\rfloor$,
that is, the set of topes
$\mathcal{K}_s{}^{\ast}-\mathcal{Q}^{\ast}$ is not a committee for
$\mathcal{M}$, a contradiction. Thus, $\mathcal{K}_s{}^{\ast}$ is
minimal. \hfill $\vartriangleleft$

Claims~1 and~2 show that the committee $\mathcal{K}^{\ast}$ is
critical.
\end{proof}

\begin{algorithm} \label{alg:1} $\quad$

\noindent {\it 01}
$\mathcal{K}_0{}^{\ast}\leftarrow\{\mathrm{T}^{(+)}\}$;

\noindent {\it 02} {\bf for $i\leftarrow 1$ to $s$}

\noindent {\it 03} \indent\indent\hspace{-2mm}{\bf do}
$\mathcal{K}_i{}^{\ast}\leftarrow\text{\tt empty set}$;

\noindent \phantom{{\it 03}}
\indent\indent\indent\hspace{-0.5mm}$\text{\em
FOUND}\leftarrow\text{\tt false}$;

\noindent {\it 04} \indent\indent\indent\hspace{-0.5mm}{\bf while}
$|\mathcal{K}_{i-1}{}^{\ast}|>0$

\noindent {\it 05} \indent\indent\indent\indent{\bf do} pick a
tope $K\in\mathcal{K}_{i-1}{}^{\ast}$;

\noindent {\it 06}
\indent\indent\indent\indent\indent\indent\hspace{-2.5mm}{\bf if}
\begin{equation} \label{eq:1}
K(j_i)=+\ \ \text{\rm and}\ \ {}_{-j_i}(-K)\in\mathcal{T}^{i-1}
\end{equation}

\noindent {\it 07}
\indent\indent\indent\indent\indent\indent\indent\hspace{-2.5mm}{\bf
then} $\text{\em FOUND}\leftarrow\text{\tt true}$;

\noindent {\it 08}
\indent\indent\indent\indent\indent\indent\indent\indent\indent\indent
\hspace{-5mm}{\bf if} there is a tope
$S\in\mathcal{K}_{i-1}{}^{\ast}$ such that
\begin{equation}
\label{eq:2}\quad\quad\quad\quad
\quad\quad\quad\quad\quad\quad\quad\quad\quad S(j_i)=+\ \
\text{\rm and}\ \ {}_{-j_i}S=-K
\end{equation}

\noindent {\it 09} \indent\indent\indent\indent
\indent\indent\indent\indent\indent\indent\hspace{-0.7mm}{\bf
then}
$\mathcal{K}_{i-1}{}^{\ast}\leftarrow\mathcal{K}_{i-1}{}^{\ast}-\{S\}$;

\noindent {\it 10} \indent\indent\indent\indent
\indent\indent\indent\indent\indent\indent\hspace{-0.7mm}{\bf
else} $\mathcal{K}_i{}^{\ast}\leftarrow\mathcal{K}_i{}^{\ast}\
\dot{\cup}\ \{K\}$;

\noindent {\it 11}
\indent\indent\indent\indent\indent\indent\indent\hspace{-2.5mm}{\bf
else} $\mathcal{K}_i{}^{\ast} \leftarrow\mathcal{K}_i{}^{\ast}\
\dot{\cup}\ \left\{{}_{-j_i}K\right\}$;

\noindent {\it 12}
\indent\indent\indent\indent\indent\indent\hspace{-2.5mm}$\mathcal{K}_{i-1}{}^{\ast}
\leftarrow\mathcal{K}_{i-1}{}^{\ast}-\{K\}$;

\noindent {\it 13} \indent\indent\indent\hspace{-0.5mm}{\bf if}
$\text{\em FOUND}=\text{\tt false}$

\noindent {\it 14} \indent\indent\indent\indent{\bf then} pick the
pair of topes $\{T'',T'''\}\subset\mathcal{T}^{i-1}$ such that
\begin{equation*}
\quad\quad\quad\quad T''(j_i)=T'''(j_i)=-\ \ \text{\rm and}\ \
{}_{-j_i}T''=-T'''\text{\em ;}
\end{equation*}

\noindent \phantom{{\it 14}}
\indent\indent\indent\indent\indent\indent\indent\hspace{-2.5mm}$\mathcal{K}_i{}^{\ast}
\leftarrow\mathcal{K}_i{}^{\ast}\ \dot{\cup}\
\left\{{}_{-j_i}T'',\ {}_{-j_i}T'''\right\}$;
\end{algorithm}

\newpage

\begin{proposition}
\label{p:1} Let $\mathcal{M}:=(E,\mathcal{T})$ be a simple
oriented matroid of rank $2$.
\begin{itemize}
\item[\rm(i)]
The set
\begin{equation}
\label{eq:5} \mathcal{K}^{\ast}:=\bigl\{K\in\mathcal{T}:\
T\in\mathcal{T},\ \mathbf{S}(K,T)=\{e\}\ \Longrightarrow\
K(e)=+\bigr\}
\end{equation}
is a critical tope committee for $\mathcal{M}$.
\item[\rm(ii)]
Committee~{\rm(\ref{eq:5})} is the set
$\mathcal{K}^{\ast}=\bmax^+(\mathcal{T})$.
\end{itemize}
\end{proposition}

\begin{proof} (i) If $\mathcal{M}$ is acyclic, then the one-element set
$\{\mathrm{T}^{(+)}\}$ is its critical tope committee.

Suppose that $\mathcal{M}$ is not acyclic. If $\mathcal{N}^0$ is
an acyclic reorientation of $\mathcal{M}$ and
$J:=(j_1,\ldots,j_s)\subset E_m$ is an ordered set of integers
such that $\mathcal{M}={}_{-J}\mathcal{N}^0$, then the proposition
follows from Lemma~\ref{p:6}.

Assertion~(ii) follows from~(i) and from Remark~\ref{p:8}(ii).
\end{proof}

\subsection{Rank $\geq 2$}

We now show that a simple oriented matroid of arbitrary rank has a
tope committee. We again use the technique of reorientations of an
initial acyclic oriented matroid $\mathcal{N}^0$. To simplify our
presentation, we will suppose below that the $i$th reorientation
of $\mathcal{N}^0$ is the oriented matroid
${}_{-[1,i]}\mathcal{N}^0$.

\begin{lemma}
\label{p:7} Let
$\mathcal{N}^0:=(E_m,\mathcal{L}^0)=(E_m,\mathcal{T}^0)$ be a
simple acyclic oriented matroid on the ground set $E_m$, with set
of covectors $\mathcal{L}^0$ and set of topes $\mathcal{T}^0$.

Let $s$ be an integer, $s\leq m$. For every $i$, $1\leq i\leq s$,
define the reorientation
$\mathcal{N}^{i}:={}_{-[1,i]}\mathcal{N}^0$.

The reorientation $\mathcal{N}^{s}$ of $\mathcal{N}^0$ has a tope
committee.
\end{lemma}

\begin{proof}
The set of covectors $\mathcal{K}_s{}^{\ast}\subset\mathcal{T}^s$
built by Algorithm~\ref{alg:3} is a tope committee for
$\mathcal{N}^s$:

{\sl Claim.} {\em The set $\mathcal{K}_s{}^{\ast}$ is of odd
cardinality, and for any $e\in E_m$, it holds
\begin{equation}
\label{eq:6} |\{K\in\mathcal{K}_s{}^{\ast}:\ K(e)=+\}|
=\left\lceil\tfrac{|\mathcal{K}_s{}^{\ast}|}{2}\right\rceil\ .
\end{equation}
}

$\triangleright$ Fix a maximal chain
$\mathbf{m}:=(R^0:=\mathrm{T}^{(+)}\ \precdot\ R^1\ \precdot\
\cdots\ \precdot\ R^m:=\mathrm{T}^{(-)})$ in the tope poset
$\mathcal{T}^0(\mathcal{L}^0,\mathrm{T}^{(+)})$. Then assign to
every tope $R^i\in\{R^1,\ldots,R^m\}$ a label
$\mathfrak{l}_i\in[1,m]$ defined by
\begin{equation}
\label{eq:14} \{\mathfrak{l}_i\}:=\mathbf{S}(R^{i-1},R^i)\ .
\end{equation}

Note that
\begin{itemize}
\item
$|\ \bigl\{R^i:\ i\in[1,m-1]\bigr\}\ \cap\
\bigl\{{}_{-\mathfrak{l}_i}(-R^i):\ i\in[2,m]\bigr\}\ |=0$ because
in the poset
$\mathcal{T}^0(\mathcal{L}^0,\mathrm{T}^{(+)})-\{\mathrm{T}^{(+)},\mathrm{T}^{(-)}\}$
its maximal chains
\begin{equation*}
\mathbf{m}-\{R^0,R^m\}
\end{equation*}
and
\begin{equation*}
\bigl(\ {}_{-\mathfrak{l}_{m}}(-R^{m})\ \precdot\
{}_{-\mathfrak{l}_{m-1}}(-R^{m-1})\ \precdot\ \cdots\ \precdot\
{}_{-\mathfrak{l}_{2}}(-R^{2})\ \bigr)
\end{equation*}
are disjoint;

\item
for every $i\in[1,m-1]$ it holds
\begin{equation*}
-R^i={}_{-\mathfrak{l}_{i+1}}(-R^{i+1})\ ;
\end{equation*}
as a consequence, the set of topes
\begin{equation*}
\bigl\{R^i:\ i\in[1,m-1]\bigr\}\ \dot\cup\
\bigl\{{}_{-\mathfrak{l}_i}(-R^i):\ i\in[2,m]\bigr\}
\end{equation*}
contains precisely $m-1$ pairs of opposites, see
Figure~\ref{fig:4};
\item
the multiset $\{R^0,{}_{-\mathfrak{l}_{1}}(-R^{1}),R^m\}=
\{\mathrm{T}^{(+)},\mathrm{T}^{(-)},\mathrm{T}^{(-)}\}$ contains exactly two pairs of opposites, namely
$\{R^0,{}_{-\mathfrak{l}_{1}}(-R^{1})\}$ and $\{R^0,R^m\}$, see Figure~\ref{fig:4};
\item
the sequence of topes
\begin{multline*}
\bigl(\mathbf{m},{}_{-\mathfrak{l}_{2}}(-R^{2}),
{}_{-\mathfrak{l}_{3}}(-R^{3}),\ldots,{}_{-\mathfrak{l}_{m}}(-R^{m}),R^0\bigr)\\=
\bigl(R^0,R^1,\ldots,R^m,{}_{-\mathfrak{l}_{2}}(-R^{2}),
{}_{-\mathfrak{l}_{3}}(-R^{3}),\ldots,{}_{-\mathfrak{l}_{m}}(-R^{m}),R^0\bigr)
\end{multline*}
is a symmetric cycle in the tope graph
$\mathcal{T}^0(\mathcal{L}^0)$ of $\mathcal{N}^0$.
\end{itemize}

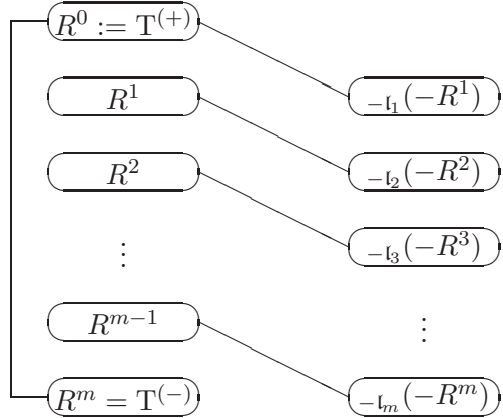
\begin{figure}[ht]
\begin{picture}(10,6)(-1.25,0.5)

\put(2,1){\oval(2,0.5)}
\put(2,1){\makebox(0,0)[c]{$R^m=\mathrm{T}^{(-)}$}}

\put(2,2){\oval(2,0.5)} \put(2,2){\makebox(0,0)[c]{$R^{m-1}$}}

\put(6,1){\oval(2,0.5)}
\put(6,1){\makebox(0,0)[c]{${}_{-\mathfrak{l}_m}(-R^m)$}}

\put(6,3){\oval(2,0.5)}
\put(6,3){\makebox(0,0)[c]{${}_{-\mathfrak{l}_3}(-R^3)$}}

\put(2,3){\makebox(0,0)[c]{$\vdots$}}
\put(6,2){\makebox(0,0)[c]{$\vdots$}}

\put(2,4){\oval(2,0.5)} \put(2,4){\makebox(0,0)[c]{$R^2$}}
\put(3,4){\line(2,-1){2}}

\put(3,2){\line(2,-1){2}}

\put(6,4){\oval(2,0.5)}
\put(6,4){\makebox(0,0)[c]{${}_{-\mathfrak{l}_2}(-R^2)$}}

\put(2,5){\oval(2,0.5)} \put(2,5){\makebox(0,0)[c]{$R^1$}}

\put(3,5){\line(2,-1){2}}

\put(6,5){\oval(2,0.5)}
\put(6,5){\makebox(0,0)[c]{${}_{-\mathfrak{l}_1}(-R^1)$}}

\put(2,6){\oval(2,0.5)}
\put(2,6){\makebox(0,0)[c]{$R^0:=\mathrm{T}^{(+)}$}}
\put(3,6){\line(2,-1){2}} \put(1,6){\line(-1,0){0.5}}
\put(0.5,6){\line(0,-1){5}} \put(0.5,1){\line(1,0){0.5}}

\end{picture}
\caption{A multiset composed of the elements $R^0:=\mathrm{T}^{(+)}\ \precdot\ R^1\ \precdot\ \cdots\
\precdot\ R^{m-1}$ of the tope poset $\mathcal{T}^0(\mathcal{L}^0,\mathrm{T}^{(+)})$ associated with a simple
acyclic oriented matroid $\mathcal{N}^0:=(E_m,\mathcal{L}^0)=(E_m,\mathcal{T}^0)$, and of their opposites.
Every pair of opposites is connected by an edge.} \label{fig:4}
\end{figure}

Now consider Algorithm~\ref{alg:4}, a slight modification of Algorithm~\ref{alg:3}. Both algorithms construct
the same committee $\mathcal{K}_s{}^{\ast}$, provided the same maximal chain of topes $\mathbf{m}$ is chosen
at their Steps~{\it 01}. Algorithm~\ref{alg:3} immediately removes pairs of opposites when they appear;
Algorithm~\ref{alg:4} removes pairs of opposites when it completes its operations at Steps~{\it 04-05}.

At Step~{\it 03}, Algorithm~\ref{alg:4} builds a multiset
$\mathcal{K}_s{}^{\ast}$ of odd cardinality that
satisfies~(\ref{eq:6}). At Steps~{\it 04-05} pairs of opposites
are thrown away; the resulting set $\mathcal{K}_s{}^{\ast}$ still
satisfies~(\ref{eq:6}). \hfill $\vartriangleleft$

\end{proof}

\newpage

\begin{algorithm} \label{alg:3} $\quad$

\noindent {\it 01}
$\mathcal{K}_0{}^{\ast}\leftarrow\{\mathrm{T}^{(+)}\}$;

\noindent \phantom{{\it 01}} $\mathbf{m}\leftarrow$ maximal chain
$(R^0\ \precdot\ R^1\ \precdot\ \cdots\ \precdot\ R^m)$ in
$\mathcal{T}^0(\mathcal{L}^0,\mathrm{T}^{(+)})$;

\noindent \phantom{{\it 01}}
$(\{\mathfrak{l}_1\},\ldots,\{\mathfrak{l}_m\})\leftarrow(\mathbf{S}(R^0,R^1),
\ldots,\mathbf{S}(R^{m-1},R^m))$;

\noindent {\it 02} {\bf for $i\leftarrow 1$ to $s$}

\noindent {\it 03} \indent\indent\hspace{-2mm}{\bf do}
$\mathcal{K}_i{}^{\ast}\leftarrow{}_{-i}(\mathcal{K}_{i-1}{}^{\ast})$;

\noindent \phantom{{\it 03}}
\indent\indent\indent\hspace{-0.5mm}$R\leftarrow R^k$ such that
$\mathfrak{l}_k=i$;

\noindent \phantom{{\it 03}}
\indent\indent\indent\hspace{-0.5mm}$\mathcal{K}_i{}^{\ast}\leftarrow
\mathcal{K}_i{}^{\ast}\ \dot\cup\ \{{}_{-[1,i]}R, {}_{-[i,m]}R\}$;

\noindent {\it 04} \indent\indent\indent\hspace{-1.6mm} {\bf
while} $\mathcal{K}_i{}^{\ast}\supset\{K,T\}$ such that $T=-K$

\noindent {\it 05} \indent\indent\indent\indent\indent {\bf do}
$\mathcal{K}_i{}^{\ast}\leftarrow\mathcal{K}_i{}^{\ast}-\{K,T\}$;
\end{algorithm}

\begin{algorithm} \label{alg:4} $\quad$

\noindent {\it 01}
$\mathcal{K}_0{}^{\ast}\leftarrow\{\mathrm{T}^{(+)}\}$;

\noindent \phantom{{\it 01}} $\mathbf{m}\leftarrow$ maximal chain
$(R^0\ \precdot\ R^1\ \precdot\ \cdots\ \precdot\ R^m)$ in
$\mathcal{T}^0(\mathcal{L}^0,\mathrm{T}^{(+)})$;

\noindent \phantom{{\it 01}}
$(\{\mathfrak{l}_1\},\ldots,\{\mathfrak{l}_m\})\leftarrow(\mathbf{S}(R^0,R^1),
\ldots,\mathbf{S}(R^{m-1},R^m))$;

\noindent {\it 02} {\bf for $i\leftarrow 1$ to $s$}

\noindent {\it 03} \indent\indent\hspace{-2mm}{\bf do} {\tt
multiset} $\mathcal{K}_i{}^{\ast}\leftarrow\{{}_{-i}K:\
K\in\mathcal{K}_{i-1}{}^{\ast}\}$;

\noindent \phantom{{\it 03}}
\indent\indent\indent\hspace{-0.5mm}$R\leftarrow R^k$ such that
$\mathfrak{l}_k=i$;

\noindent \phantom{{\it 03}}
\indent\indent\indent\hspace{-0.5mm}$\mathcal{K}_i{}^{\ast}\leftarrow
\{\mathcal{K}_i{}^{\ast}, {}_{-[1,i]}R, {}_{-[i,m]}R\}$;

\noindent {\it 04} {\bf while}
$\mathcal{K}_s{}^{\ast}\supset\{K,T\}$ such that $T=-K$

\noindent {\it 05} \indent\indent\hspace{-2.7mm} {\bf do}
$\mathcal{K}_s{}^{\ast}\leftarrow\mathcal{K}_s{}^{\ast}-\{K,T\}$;
\end{algorithm}

Proposition~\ref{p:3} of Section~\ref{app:1} asserts that
Algorithm~\ref{alg:3} in fact constructs critical committees.

\begin{example}
Consider the acyclic oriented matroid
$\mathcal{N}^0:=(E_6,\mathcal{T}^0)$ given by the central
hyperplane arrangement of Figure~{\rm\ref{fig:7}}.

If the maximal chain {\scriptsize
\begin{equation*}
\text{\normalsize$\mathbf{m}$}:=\quad
\begin{matrix}
(&R^0&:=&+&+&+&+&+&+&\\ &R^1&:=&+&+&-&+&+&+&\\
&R^2&:=&-&+&-&+&+&+&\\ &R^3&:=&-&+&-&-&+&+&\\
&R^4&:=&-&+&-&-&+&-&\\ &R^5&:=&-&-&-&-&+&-&\\
&R^6&:=&-&-&-&-&-&-&)
\end{matrix}
\end{equation*}
} in its tope poset
$\mathcal{T}^0(\mathcal{L}^0,\mathrm{T}^{(+)})$ is chosen at
Step~{\it 01} of Algorithm~{\rm\ref{alg:3}}, with
\begin{equation*}
\mathfrak{l}_1=3,\ \ \mathfrak{l}_2=1,\ \ \mathfrak{l}_3=4,\ \
\mathfrak{l}_4=6,\ \ \mathfrak{l}_5=2,\ \ \mathfrak{l}_6=5,
\end{equation*}
then the algorithm produces, when applying to the reorientation $\mathcal{N}^6:={}_{-[1,6]}\mathcal{N}^0$,
the following sequence of critical tope committees: {\scriptsize
\begin{gather*}
\mathcal{K}_1{}^{\ast}=\quad
\begin{matrix}
\{&-&+&+&+&+&+&\\ &+&+&-&+&+&+&\\ &+&-&+&-&-&-&\}
\end{matrix}\quad\quad
\mathcal{K}_2{}^{\ast}=\quad
\begin{matrix}
\{&-&-&+&+&+&+&\\ &+&-&-&+&+&+&\\ &+&+&+&-&-&-&\\ &+&+&-&-&+&-&\\
&-&+&+&+&-&+&\}
\end{matrix}\\
\mathcal{K}_3{}^{\ast}=\quad
\begin{matrix}
\{&+&-&+&+&+&+&\\ &+&+&+&-&+&-&\\ &-&+&-&+&-&+&\}
\end{matrix}\quad\quad
\mathcal{K}_4{}^{\ast}=\quad
\begin{matrix}
\{&+&+&+&+&+&-&\\ &-&+&-&-&-&+&\\ &+&-&+&+&+&+&\}
\end{matrix}\\
\mathcal{K}_5{}^{\ast}=\quad
\begin{matrix}
\{&-&+&-&-&+&+&\\ &+&-&+&+&-&+&\\ &+&+&+&+&+&-&\}
\end{matrix}\quad\quad
\mathcal{K}_6{}^{\ast}=\quad
\begin{matrix}
\{&+&+&+&+&+&+&\}
\end{matrix}
\end{gather*}
}

\end{example}

Algorithm~{\rm\ref{alg:3}} does not necessarily construct a tope
committee of minimal cardinality. For example, it builds for the
reorientation $\mathcal{N}^5:={}_{-[1,5]}\mathcal{N}^0$ of the
oriented matroid $\mathcal{N}^0$ {\rm(}which is realized by the
hyperplane arrangement of Figure~{\rm\ref{fig:7})} a three-element
committee, while $\mathcal{N}^5$ is acyclic, that is, the
one-element set $\{\mathrm{T}^{(+)}\}$ is a committee for
$\mathcal{N}^5$.

Any simple oriented matroid has a tope committee of cardinality
that is less than or equal to the cardinality of its ground set:

\begin{theorem}
Let $\mathcal{N}^0:=(E_m,\mathcal{T}^0)$ be a simple acyclic
oriented matroid on the ground set $E_m$, with set of topes
$\mathcal{T}^0$. Let $s$ be an integer, $s\leq m$. Define the
reorientation $\mathcal{N}^{s}:={}_{-[1,s]}\mathcal{N}^0$.

\begin{itemize}
\item[\rm(i)] $\mathcal{N}^{s}$ has a
tope committee $\mathcal{K}_s{}^{\ast}$ such that
$|\mathcal{K}_s{}^{\ast}|\leq m$ when $m$ is odd, and
$|\mathcal{K}_s{}^{\ast}|\leq m-1$ when $m$ is even.

As a consequence, any simple oriented matroid $\mathcal{M}$ on the
set $E_m$ has a tope committee $\mathcal{K}^{\ast}$ with
\begin{equation*}
|\mathcal{K}^{\ast}|\leq\begin{cases}m\ ,&\text{if\/ $m$ is
odd,}\\ m-1\ ,&\text{if\/ $m$ is even.}
\end{cases}
\end{equation*}

\item[\rm(ii)]
If $s\leq\lfloor m/2\rfloor$ then $\mathcal{N}^{s}$ has a tope
committee $\mathcal{K}_s{}^{\ast}$ with
$|\mathcal{K}_s{}^{\ast}|\leq 1+2s$.
\end{itemize}
\end{theorem}

\begin{proof}
Apply Algorithm~\ref{alg:3} to $\mathcal{N}^s$; this algorithm constructs a tope committee for
$\mathcal{N}^s$, see the proof of Lemma~\ref{p:7}.
\begin{itemize}
\item[\rm(i)]
\begin{itemize}
\item[$\bullet$]
If $m$ is odd, then the tope committee of maximal cardinality
which can be constructed by Algorithm~\ref{alg:3} is the set
\begin{equation}
\label{eq:11} \mathcal{K}_s{}^{\ast}=\{{}_{-[1,s]}R^0\}\ \
\dot\cup\ \ \{{}_{-[1,s]}R^k,\ {}_{-[s,m]}R^k:\ 2\leq k\leq m-1,\
\text{$k$ even}\}
\end{equation}
of cardinality $m$.

Figure~\ref{fig:5} depicts such an arrangement of topes,
cf.~Figure~\ref{fig:4}.

\begin{figure}[ht]
\begin{picture}(10,6)(-1,0.5)

\put(2,1){\oval(2,0.5)} \put(2,1){\makebox(0,0)[c]{$\times$}}

\put(2,2){\oval(2,0.5)}
\put(2,2){\makebox(0,0)[c]{${}_{-[1,s]}R^{m-1}$}}

\put(6,1){\oval(2,0.5)} \put(6,1){\makebox(0,0)[c]{$\times$}}

\put(6,3){\oval(2,0.5)} \put(6,3){\makebox(0,0)[c]{$\times$}}

\put(2,3){\makebox(0,0)[c]{$\vdots$}}
\put(6,2){\makebox(0,0)[c]{$\vdots$}}

\put(2,4){\oval(2,0.5)}
\put(2,4){\makebox(0,0)[c]{${}_{-[1,s]}R^2$}}

\put(6,4){\oval(2,0.5)}
\put(6,4){\makebox(0,0)[c]{${}_{-[s,m]}R^2$}}

\put(2,5){\oval(2,0.5)} \put(2,5){\makebox(0,0)[c]{$\times$}}

\put(6,5){\oval(2,0.5)} \put(6,5){\makebox(0,0)[c]{$\times$}}

\put(2,6){\oval(2,0.5)}
\put(2,6){\makebox(0,0)[c]{${}_{-[1,s]}R^0$}}
\end{picture}
\caption{The tope committee of maximal cardinality~(\ref{eq:11})
which can be constructed by Algorithm~\ref{alg:3} in the case of
$m$ odd; cf.~Figure~\ref{fig:4}.} \label{fig:5}
\end{figure}

\item[$\bullet$]
If $m$ is even, then the tope committee of maximal cardinality
which can be constructed by Algorithm~\ref{alg:3} is either the
set
\begin{equation}
\label{eq:12} \mathcal{K}_s{}^{\ast}=\{{}_{-[1,s]}R^k,\
{}_{-[s,m]}R^k:\ 2\leq k\leq m-2,\ \text{$k$ even}\}\ \ \dot\cup\
\ \{{}_{-[s,m]}R^m\}
\end{equation}
or the set
\begin{equation}
\label{eq:13} \mathcal{K}_s{}^{\ast}=\{{}_{-[1,s]}R^1\}\ \
\dot\cup\ \ \{{}_{-[1,s]}R^k,\ {}_{-[s,m]}R^k:\ 3\leq k\leq m-1,\
\text{$k$ odd}\}\ ,
\end{equation}
see Figures~\ref{fig:6}(a) and~(b), respectively,
cf.~Figure~\ref{fig:4}. We have $|\mathcal{K}_s{}^{\ast}|=m-1$.

\begin{figure}[ht]
\begin{picture}(10,6)(-0.75,0.25)

\put(0,1){\oval(2,0.5)} \put(0,1){\makebox(0,0)[c]{$\times$}}

\put(0,2){\oval(2,0.5)} \put(0,2){\makebox(0,0)[c]{$\times$}}

\put(2.5,1){\oval(2,0.5)}
\put(2.5,1){\makebox(0,0)[c]{${}_{-[s,m]}R^m$}}

\put(2.5,3){\oval(2,0.5)} \put(2.5,3){\makebox(0,0)[c]{$\times$}}

\put(0,3){\makebox(0,0)[c]{$\vdots$}}
\put(2.5,2){\makebox(0,0)[c]{$\vdots$}}

\put(0,4){\oval(2,0.5)}
\put(0,4){\makebox(0,0)[c]{${}_{-[1,s]}R^2$}}

\put(2.5,4){\oval(2,0.5)}
\put(2.5,4){\makebox(0,0)[c]{${}_{-[s,m]}R^2$}}

\put(0,5){\oval(2,0.5)} \put(0,5){\makebox(0,0)[c]{$\times$}}

\put(2.5,5){\oval(2,0.5)} \put(2.5,5){\makebox(0,0)[c]{$\times$}}

\put(0,6){\oval(2,0.5)} \put(0,6){\makebox(0,0)[c]{$\times$}}
\put(1.25,0.5){\makebox(0,0)[t]{\bf (a)}}

\put(4.25,1){\line(0,1){5}}

\put(6,1){\oval(2,0.5)} \put(6,1){\makebox(0,0)[c]{$\times$}}

\put(6,2){\oval(2,0.5)}
\put(6,2){\makebox(0,0)[c]{${}_{-[1,s]}R^{m-1}$}}

\put(8.5,1){\oval(2,0.5)} \put(8.5,1){\makebox(0,0)[c]{$\times$}}

\put(8.5,3){\oval(2,0.5)}
\put(8.5,3){\makebox(0,0)[c]{${}_{-[s,m]}R^3$}}

\put(6,3){\makebox(0,0)[c]{$\vdots$}}
\put(8.5,2){\makebox(0,0)[c]{$\vdots$}}

\put(6,4){\oval(2,0.5)} \put(6,4){\makebox(0,0)[c]{$\times$}}

\put(8.5,4){\oval(2,0.5)} \put(8.5,4){\makebox(0,0)[c]{$\times$}}

\put(6,5){\oval(2,0.5)}
\put(6,5){\makebox(0,0)[c]{${}_{-[1,s]}R^1$}}

\put(8.5,5){\oval(2,0.5)} \put(8.5,5){\makebox(0,0)[c]{$\times$}}

\put(6,6){\oval(2,0.5)} \put(6,6){\makebox(0,0)[c]{$\times$}}
\put(7.25,0.5){\makebox(0,0)[t]{\bf (b)}}
\end{picture}
\caption{The tope committees of maximal cardinality which can be
constructed by Algorithm~\ref{alg:3} in the case of $m$ even:
(\ref{eq:12}) on the left, and (\ref{eq:13}) on the right;
cf.~Figure~\ref{fig:4}.} \label{fig:6}
\end{figure}
\end{itemize}

Figure~\ref{fig:4} suggests that if Algorithm~\ref{alg:3} builds
for $s=\lfloor m/2\rfloor$ a tope committee of the
form~(\ref{eq:11}),~(\ref{eq:12}) or~(\ref{eq:13}), then for every
$s$ such that $s>\lfloor m/2\rfloor$, the cardinality of
$\mathcal{K}_s{}^{\ast}$ will decrease because pairs of opposites
will be removed.

\item[\rm(ii)] The committee constructions of maximal cardinality
which we have considered in part~(i) were built under $s=\lfloor
m/2\rfloor$. One can argue in an analogous way to prove
assertion~(ii).
\end{itemize}
\end{proof}

\section{Graphs Related to Tope Committees}
\label{section:4}

From the point of view of modeling of decision-making procedures, naturally associated with a simple oriented
matroid $\mathcal{M}:=(E_m,\mathcal{T})$ is a certain graph $\boldsymbol{\Gamma}(\mathcal{M})$ that is
isomorphic to the Kneser graph $\KG\bigl(\{T^-:\ T\in\mathcal{T}\}\bigr)$ of the family of the negative parts
of the topes of $\mathcal{M}$; the vertex sets of odd cycles in $\boldsymbol{\Gamma}(\mathcal{M})$ are tope
committees:

\begin{lemma}
Let $\mathcal{M}:=(E_m,\mathcal{L})=(E_m,\mathcal{T})$ be a simple
oriented matroid.

Consider the graph
$\boldsymbol{\Gamma}:=\boldsymbol{\Gamma}(\mathcal{M})$ defined by
\begin{equation}
\label{eq:21}
\begin{split}
\mathfrak{V}(\boldsymbol{\Gamma})\ \ &\ :=\ \ \mathcal{T}\ ,\\
\{T^k,T^l\}\in\mathfrak{E}(\boldsymbol{\Gamma})\ \
&\Longleftrightarrow\ \ (T^k)^+\cup(T^l)^+=E_m\ .
\end{split}
\end{equation}

If\/ $\pmb{C}$ is an odd cycle in $\boldsymbol{\Gamma}$, then the
set of its vertices $\mathfrak{V}(\pmb{C})$ is a tope committee
for $\mathcal{M}$.
\end{lemma}

\begin{proof}
Assume that there is an element $e\in E_m$ such that
$|\{K\in\mathfrak{V}(\pmb{C}):\
K(e)=-\}|\geq\left\lceil\tfrac{|\mathfrak{V}(\pmb{C})|}{2}\right\rceil$.
Then there exists an edge $\{T',T''\}\in\mathfrak{E}(\pmb{C})$
with $(T')^+\cup(T'')^+\not\ni e$; hence $(T')^+\cup(T'')^+\neq
E_m$, a contradiction.
\end{proof}

\begin{example} Consider the reorientation $\mathcal{N}^3:=
{}_{-[1,3]}\mathcal{N}^0$ of the oriented matroid $\mathcal{N}^0$
which is realized by the hyperplane arrangement of
Figure~{\rm\ref{fig:7}}. The set of vertices of the $5$-cycle in
$\boldsymbol{\Gamma}(\mathcal{N}^3)$, shown in
Figure~{\rm\ref{fig:10}}, is a tope committee for $\mathcal{N}^3$.

\begin{figure}[ht]
\begin{picture}(12,5)(0,0)
\put(1.5,2.5){\oval(2,0.5)} \put(5,0.5){\oval(2,0.5)}
\put(5,4.5){\oval(2,0.5)} \put(10,1){\oval(2,0.5)}
\put(10,4){\oval(2,0.5)}

\put(1.5,2.5){\makebox(0,0)[c]{\scriptsize $-+-+-+$}}
\put(5,0.5){\makebox(0,0)[c]{\scriptsize $+-++++$}}
\put(5,4.5){\makebox(0,0)[c]{\scriptsize $+++-+-$}}
\put(10,1){\makebox(0,0)[c]{\scriptsize $+++---$}}
\put(10,4){\makebox(0,0)[c]{\scriptsize $--++++$}}

\put(2.45,2.7){\line(1,1){1.6}} \put(5.95,4.7){\line(5,-1){3}}
\put(10.5,3.75){\line(0,-1){2.5}} \put(5.95,0.3){\line(5,1){3}}
\put(2.45,2.3){\line(1,-1){1.6}}
\end{picture}
\caption{A $5$-cycle in the graph $\boldsymbol{\Gamma}(\mathcal{N}^3)$, defined by~{\rm(\ref{eq:21})}, that
is associated with the reorientation $\mathcal{N}^3:={}_{-[1,3]}\mathcal{N}^0$. The oriented matroid
$\mathcal{N}^0$ is realized by the hyperplane arrangement of Figure~\ref{fig:7}. The set of vertices of the
cycle is a tope committee for $\mathcal{N}^3$.} \label{fig:10}
\end{figure}
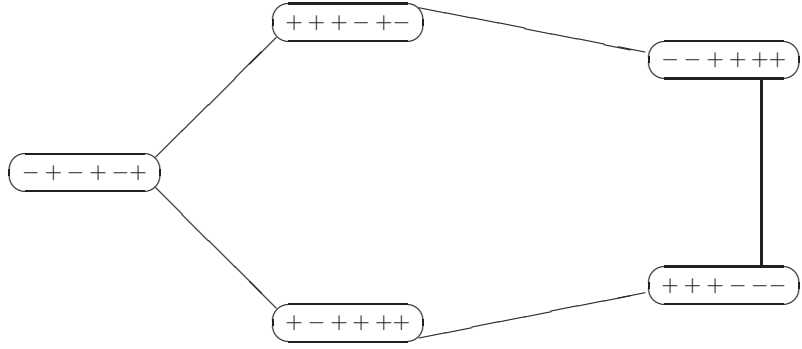
\end{example}

\subsection{Symmetric Cycles in the Tope Graph} \label{app:1}

We now show that a direct graph generalization of centrally-symmetric cycles of adjacent regions in
hyperplane arrangements from~\cite{AK-IEEE} leads to constructions of odd cycles in subgraphs of the graphs
$\boldsymbol{\Gamma}$ defined by (\ref{eq:21}).

\begin{proposition}
\label{p:2} Let $\mathcal{M}:=(E_m,\mathcal{L})=(E_m,\mathcal{T})$
be a simple oriented matroid. Let
$\pmb{R}:=(T^0,T^1,\ldots,T^{2m-1},T^0)$ be a symmetric cycle
{\rm(}that does not contain the positive tope
$\mathrm{T}^{(+)}${\rm)} in the tope graph
$\mathcal{T}(\mathcal{L})$.

Consider the graph $\pmb{G}$ defined by
\begin{equation}
\label{eq:8}
\begin{split}
\mathfrak{V}(\pmb{G})\ \ &\ :=\ \
\mathfrak{V}(\pmb{R})=\{T^0,T^1,\ldots,T^{2m-1}\}\ ,\\
\{T^k,T^l\}\in\mathfrak{E}(\pmb{G})\ \ &\Longleftrightarrow\ \
(T^k)^+\cup(T^l)^+=E_m\ .
\end{split}
\end{equation}

The set $\bmax^+\bigl(\mathfrak{V}(\pmb{R})\bigr)$ is the vertex
set of an odd cycle in $\pmb{G}$.
\end{proposition}

\begin{proof}
We without loss of generality suppose that
$T^0\in\bmax^+\bigl(\mathfrak{V}(\pmb{R})\bigr)$.

Example~\ref{e:777} and Figure~\ref{fig:8} illustrate the proof.

Note that the path $(T^1,T^2,\ldots,T^m)$ contains at least one vertex $T^j$ such that
$T^j\in\bmax^+\bigl(\mathfrak{V}(\pmb{R})\bigr)$. This follows from the observation that $|(T^0)^+|$
$<m=|\mathbf{S}(T^0,T^m)|$ and $|(T^m)^+|>0$, and from Remark~\ref{p:8}(ii).

Let $T^l$ be a vertex of $\pmb{R}$ such that $1<l<m$,
$T^l\in\bmax^+\bigl(\mathfrak{V}(\pmb{R})\bigr)$, and
$(T^l)^+\supsetneqq(T^m)^+=(-T^0)^+$. The pair $\{T^0,T^l\}$ is an
edge of $\pmb{G}$.

We have $\{T^0,T^j\}\in\mathfrak{E}(\pmb{G})$ for all $j$, $l\leq
j\leq m$.

On the contrary, if $0<j<l$ then
$\{T^0,T^j\}\not\in\mathfrak{E}(\pmb{G})$. Indeed, let $\{e\}$ be
the one-element separation set for the topes $T^l$ and $T^{l-1}$.
Then we have $e\not\in(T^j)^+$ and $e\not\in(T^0)^+$.

Similarly, there is a unique vertex $T^p$ of the cycle $\pmb{R}$
such that $p>m$, $T^p\in\bmax^+\bigl(\mathfrak{V}(\pmb{R})\bigr)$,
and $\{T^0,T^p\}\in\mathfrak{E}(\pmb{G})$. Here
$(T^p)^+\supsetneqq(T^m)^+=(-T^0)^+$. For all $j$, $m\leq j\leq
p$, we have $\{T^0,T^j\}\in\mathfrak{E}(\pmb{G})$. If $p<j\leq
2m-1$ then $\{T^0,T^j\}\not\in\mathfrak{E}(\pmb{G})$.

Thus,
\begin{multline*}
\bigl\{\{T', T''\}\in\mathfrak{E}(\pmb{G}):\ \{T',T''\}\ni
T^0\bigr\}\\=\bigl\{\{T^0,T^l\},\{T^0,T^{l+1}\}\ldots,
\{T^0,T^{m-1}\},\{T^0,T^m\},\\
\{T^0,T^{m+1}\},\ldots,\{T^0,T^{p-1}\},\{T^0,T^p\}\bigr\}\ .
\end{multline*}

Note that for all $j$, $l<j<p$, we have $T^j\not\in
\bmax^+\bigl(\mathfrak{V}(\pmb{R})\bigr)$.

Let $T^i$ be a vertex of $\pmb{R}$ such that $1<i\leq l$, $T^i\in
\bmax^+\bigl(\mathfrak{V}(\pmb{R})\bigr)$, and
$(T^{i-1})^+\subsetneqq(T^0)^+$. The pair $\{T^p,T^i\}$ is an edge
of $\pmb{G}$.

We have the inclusion
\begin{equation*}
\{T^0,T^l\},\ \{T^0,T^p\},\ \{T^p,T^i\}\ \in \
\mathfrak{E}(\pmb{G})\ .
\end{equation*}

If $i=l$ then the sequence of vertices $(T^0,T^l,T^p,T^0)$ is a
triangle in $\pmb{G}$ with the property
$\{T^0,T^l,T^p\}=\bmax^+\bigl(\mathfrak{V}(\pmb{R})\bigr)$.

In the general case, when $i\leq l$, consider successively all the
vertices $T^j$ with $i\leq j\leq l$ to see that the set
$\bmax^+\bigl(\mathfrak{V}(\pmb{R})\bigr)$ is of odd cardinality
because
\begin{multline}
\label{eq:9} \bigl|\ \{i:\ 1<i<m,\
T^i\in\bmax^+\bigl(\mathfrak{V}(\pmb{R})\bigr)\ \bigr|\\=\bigl|\
\{i:\ m<i<2m-1,\ T^i\in\bmax^+\bigl(\mathfrak{V}(\pmb{R})\bigr)\
\bigr|\ ,
\end{multline}
and this set is the vertex set $\mathfrak{V}(\pmb{C})$ of a cycle
$\pmb{C}$ in $\pmb{G}$: if
\begin{equation*}
\bmax^+\bigl(\mathfrak{V}(\pmb{R})\bigr)=\{T^0, T^{k_1},\ldots,T^{k_d},T^{k_{d+1}},\ldots,T^{k_{2d}}\}\ ,\ \
\ 0<k_1<\cdots<k_{2d}\ ,
\end{equation*}
then the family of edges of this cycle is
\begin{equation}
\label{eq:10}
\begin{split}
\mathfrak{E}(\pmb{C})\ =\ \big\{\ \{T^0,T^{k_d}\}\ &,\ \
\{T^0,T^{k_{d+1}}\}\ ,\\ \{T^{k_1},T^{k_{d+1}}\}\ &,\ \
\{T^{k_1},T^{k_{d+2}}\}\ ,\ \ \ldots\ ,\\
\{T^{k_{d-1}},T^{k_{2d-1}}\}\ &,\ \ \{T^{k_{d-1}},T^{k_{2d}}\}\
,\\ \quad\ &\phantom{,}\ \ \{T^{k_d},T^{k_{2d}}\}\ \big\}\ .
\end{split}
\end{equation}
\end{proof}

\begin{example}
\label{e:777} Consider the reorientation $\mathcal{M}:={}_{-[1,2]}\mathcal{N}^0$ of an oriented matroid
$\mathcal{N}^0$, where $\mathcal{N}^0$ is realized by the hyperplane arrangement of Figure~{\rm\ref{fig:7}}.
Figure~{\rm\ref{fig:8}} depicts the set of vertices of a symmetric cycle $\pmb{R}$ in the tope graph of
$\mathcal{M}$, and the corresponding graph $\pmb{G}$ defined by~{\rm(\ref{eq:8})}. The set
$\bmax^+\bigl(\mathfrak{V}(\pmb{R})\bigr)$ is the vertex set of the odd cycle in $\pmb{G}$.

\begin{figure}[ht]
\begin{picture}(12,5)(0,0)
\put(1.5,2){\oval(2,0.5)} \put(1.5,3){\oval(2,0.5)}
\put(2.5,1){\oval(2,0.5)} \put(2.5,4){\oval(2,0.5)}
\put(5,0.5){\oval(2,0.5)} \put(5,4.5){\oval(2,0.5)}
\put(7.5,0.5){\oval(2,0.5)} \put(7.5,4.5){\oval(2,0.5)}
\put(10,1){\oval(2,0.5)} \put(10,4){\oval(2,0.5)}
\put(11,2){\oval(2,0.5)} \put(11,3){\oval(2,0.5)}

\put(1.5,2){\makebox(0,0)[c]{\scriptsize $--++-+$}}

\put(1,1.5){\makebox(0,0)[c]{\small $T^1$}}

\put(1.5,3){\makebox(0,0)[c]{\scriptsize $--++++$}}

\put(1,3.5){\makebox(0,0)[c]{\small $T^i$}}

\put(2.5,1){\makebox(0,0)[c]{\scriptsize $-+++-+$}}

\put(2,0.5){\makebox(0,0)[c]{\small $T^0$}}

\put(2.5,4){\makebox(0,0)[c]{\scriptsize $---+++$}}
\put(5,0.5){\makebox(0,0)[c]{\scriptsize $-+++--$}}

\put(4.75,0){\makebox(0,0)[c]{\small $T^{2m-1}$}}

\put(5,4.5){\makebox(0,0)[c]{\scriptsize $+--+++$}}

\put(5.5,5){\makebox(0,0)[c]{\small $T^l$}}

\put(7.5,0.5){\makebox(0,0)[c]{\scriptsize $-++---$}}
\put(7.5,4.5){\makebox(0,0)[c]{\scriptsize $+---++$}}
\put(10,1){\makebox(0,0)[c]{\scriptsize $+++---$}}
\put(10,4){\makebox(0,0)[c]{\scriptsize $+---+-$}}

\put(10.5,4.5){\makebox(0,0)[c]{\small $T^m$}}

\put(11,2){\makebox(0,0)[c]{\scriptsize $++----$}}
\put(11,3){\makebox(0,0)[c]{\scriptsize $++--+-$}}

\put(11.5,3.5){\makebox(0,0)[c]{\small $T^p$}}

\put(2.5,3){\line(1,0){7.5}} \put(2.4,2.75){\line(4,-1){6.6}}
\put(3.45,1.2){\line(2,5){1.22}} \put(3.5,1){\line(4,1){6.9}}

\put(3.5,1.1){\line(1,1){3.15}}

\put(5,4.22){\line(3,-2){4.42}} \put(4.8,0.75){\line(0,1){3.5}}
\end{picture}
\caption{The graph $\pmb{G}$, defined by~{\rm(\ref{eq:8})}, that
corresponds to a symmetric cycle
$\pmb{R}:=(T^0,T^1,\ldots,T^{2m-1},T^0)$ in the tope graph of the
reorientation $\mathcal{M}:={}_{-[1,2]}\mathcal{N}^0$, where the
oriented matroid $\mathcal{N}^0$ is realized by the hyperplane
arrangement of Figure~\ref{fig:7}; $m=6$. The edges that connect
opposites are {\bf not} depicted. The set
$\bmax^+\bigl(\mathfrak{V}(\pmb{R})\bigr)$ is the vertex set of
the $5$-cycle in $\pmb{G}$.} \label{fig:8}
\end{figure}
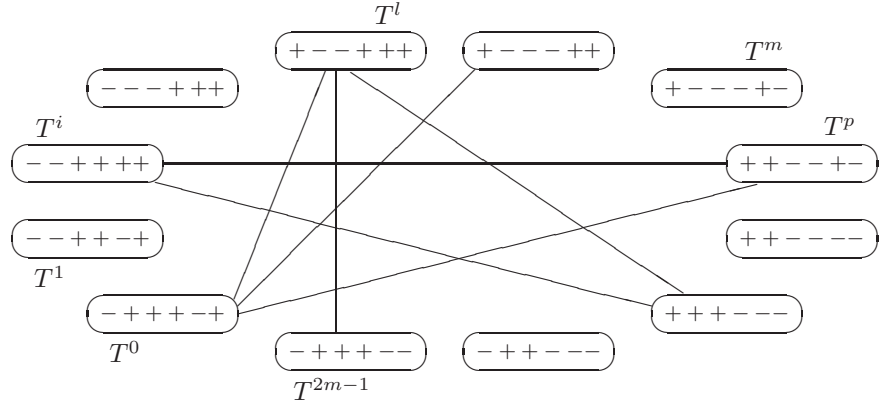

\end{example}

\begin{lemma}
\label{p:5} Let $\mathcal{M}:=(E_m,\mathcal{L})=(E_m,\mathcal{T})$
be a simple oriented matroid. Let
$\pmb{R}:=(T^0,T^1,\ldots,T^{2m-1},T^0)$ be a symmetric cycle
{\rm(}that does not contain the positive tope
$\mathrm{T}^{(+)}${\rm)} in the tope graph
$\mathcal{T}(\mathcal{L})$. For every $e\in E_m$, we have
\begin{equation*}
\bigl|\ \bigl\{T\in\bmax^+\bigl(\mathfrak{V}(\pmb{R})\bigr):\
T(e)=+\bigl\}\ \bigr|=\left\lceil\frac{|\
\bmax^+\bigl(\mathfrak{V}(\pmb{R})\bigr)\ |}{2}\right\rceil\ .
\end{equation*}
\end{lemma}

\begin{proof}
Let $\pmb{G}$ and $\pmb{C}$ be the graph and the odd cycle,
respectively, which were constructed in the proof of
Proposition~\ref{p:2}; see descriptions~(\ref{eq:8})
and~(\ref{eq:10}).

Let $\pmb{P}:=(T^{l_0},\ldots,T^{l_{m-1}})$ be the $(m-1)$-path in $\pmb{R}$ such that $T^{l_0}(e)$
$=\cdots=T^{l_{m-1}}(e)=+$; see Remark~\ref{p:8}(i).

Without loss of generality suppose that $(T^0)^+\supseteq(T^{l_0})^+\supsetneqq(T^{k_d})^-$. Then we have
$T^{k_d}\in\mathfrak{V}(\pmb{P})$, and the assertion follows from~(\ref{eq:9}) because $\bigl\{T$
$\in\bmax^+\bigl(\mathfrak{V}(\pmb{R})\bigr):\
T(e)=+\bigl\}=\mathfrak{V}(\pmb{P})\cap\bmax^+\bigl(\mathfrak{V}(\pmb{R})\bigr)=
\{T^0,T^{k_1},T^{k_2},\ldots,$ $T^{k_d}\}$.
\end{proof}

The higher-rank analogue of Proposition~\ref{p:1} is as follows:

\begin{proposition}
\label{p:4} Let $\mathcal{M}:=(E_m,\mathcal{L})=(E_m,\mathcal{T})$
be a simple oriented matroid. Let
$\pmb{R}:=(T^0,T^1,\ldots,T^{2m-1},T^0)$ be a symmetric cycle in
the tope graph $\mathcal{T}(\mathcal{L})$. The set
\begin{equation}
\label{eq:7} \mathcal{K}^{\ast}:=
\bmax^+\bigl(\mathfrak{V}(\pmb{R})\bigr)
\end{equation}
or, equivalently,
\begin{equation}
\label{eq:20} \mathcal{K}^{\ast}=\bigl\{\
T\in\mathfrak{V}(\pmb{R}):\ S\in\mathfrak{V}(\pmb{R})\ ,\
\mathbf{S}(S,T)=\{e\}\ \Longrightarrow\ T(e)=+\ \bigr\}
\end{equation}
is a critical tope committee for $\mathcal{M}$, that satisfies the
equality
\begin{equation}
\label{eq:36} |\{K\in\mathcal{K}{}^{\ast}:\ K(e)=+\}|
=\left\lceil\tfrac{|\mathcal{K}{}^{\ast}|}{2}\right\rceil\ ,
\end{equation}
for every $e\in E$.
\end{proposition}

\begin{proof} Descriptions~(\ref{eq:7}) and~(\ref{eq:20}) are
equivalent by Remark~\ref{p:8}(ii).

If $\mathrm{T}^{(+)}\in\mathfrak{V}(\pmb{R})$ then~(\ref{eq:7}) is
the one-element set $\{\mathrm{T}^{(+)}\}$, that is a critical
committee for $\mathcal{M}$; we are done.

If $\mathrm{T}^{(+)}\not\in\mathfrak{V}(\pmb{R})$ then
Lemma~\ref{p:5} implies that $\mathcal{K}^{\ast}$ is a tope
committee that satisfies~(\ref{eq:36}), for all $e\in E_m$. We
have to show that $\mathcal{K}^{\ast}$ is critical.

Assume that there is a proper subset $\mathcal{Q}^{\ast}$ of the
set $\mathcal{K}^{\ast}$ such that
$\mathcal{K}^{\ast}-\mathcal{Q}^{\ast}$ is a committee for
$\mathcal{M}$. Since
$\mathrm{T}^{(+)}\not\in\mathfrak{V}(\pmb{R})$, we have
$|\mathcal{K}^{\ast}-\mathcal{Q}^{\ast}|>1$.

We without loss of generality suppose that $T^0\in\bmax^+\bigl(\mathfrak{V}(\pmb{R})\bigr)$ and $T^0$
$\not\in\mathcal{Q}^{\ast}$.

Let $\{g\}:=\mathbf{S}(T^0,T^1)$ and
$\{f\}:=\mathbf{S}(T^{2m-1},T^0)$; note that $f\neq g$. We have
$T^0(f)=T^0(g)=+$. For every $k$, $0<k<m$, we have $T^k(f)=+$; for
every $l$, $m<l<2m$, we have $T^l(g)=+$, see Remark~\ref{p:8}(i).

By Lemma~\ref{p:5}, for every $e\in\{f,g\}$ it holds $|\{Q\in\mathcal{Q}^{\ast}:\ Q(e)=-\}|$
$=|\{R\in\mathcal{Q}^{\ast}:\ R(e)=+\}|$, and we have
\begin{equation*}
|\{T^1,T^2,\ldots,T^{m-1}\}\cap\mathcal{Q}^{\ast}|=
|\{T^{m+1},T^{m+2},\ldots,T^{2m-1}\}\cap\mathcal{Q}^{\ast}|=
\frac{|\mathcal{Q}^{\ast}|}{2}\ ;
\end{equation*}
$\mathcal{Q}^{\ast}$ is of even cardinality.

Let $T^{k_2}$ be a tope from the set $\mathcal{Q}^{\ast}$, and
$(T^{k_1},T^{k_2},T^{k_3})$ a $2$-path in $\pmb{R}$. If
$\{p\}:=\mathbf{S}(T^{k_1},T^{k_2})$ and
$\{q\}:=\mathbf{S}(T^{k_2},T^{k_3})$, then there is an element
$h\in\{p,q\}$ such that
$|\{K\in\mathcal{K}^{\ast}-\mathcal{Q}^{\ast}:\
K(h)=+\}|=\left\lfloor\tfrac{|\mathcal{K}^{\ast}-\mathcal{Q}^{\ast}|}{2}
\right\rfloor$,
that is, the set of topes $\mathcal{K}^{\ast}-\mathcal{Q}^{\ast}$
is not a committee for $\mathcal{M}$, a contradiction. Thus,
$\mathcal{K}^{\ast}$ is minimal and, as a consequence, it is
critical, in view of~(\ref{eq:36}).
\end{proof}

Note that under the hypothesis of the Proposition~\ref{p:4} the
set $-\bigl(\
\mathfrak{V}(\pmb{R})-\bmax^+\bigl(\mathfrak{V}(\pmb{R})\bigr)\
\bigr)$ is a tope committee for $\mathcal{M}$ as well.

We now discuss some poset-theoretic properties of topes which are
useful for analysis of the coverings of the ground sets of
oriented matroids by pairs of the positive parts of topes.

\begin{corollary}\label{prop:36} Let
$\mathcal{M}:=(E_m,\mathcal{L})=(E_m,\mathcal{T})$ be a simple
oriented matroid that is not acyclic. Let $\mathbf{m}$ be an
arbitrary maximal chain in the tope poset
$\mathcal{T}(\mathcal{L},B)$ with base tope
$B\in\bmax^+(\mathcal{T})$.

\begin{itemize}
\item[\rm(i)]
Let $c:=\max\{|T^+|:\ T\in\mathcal{T}\}$.

The subchain $\bmax^+(\mathbf{m})$ contains a unique tope $K$ such
that $B^+\cup K^+=E_m$. The poset rank $\rho(K)$ of $K$ satisfies
the inequality
\begin{equation}
\label{eq:31} \rho(K)\geq 2m-c-|B^+|\ .
\end{equation}

For a tope $R\in\mathbf{m}$, we have
\begin{equation}
\label{eq:30} B^+\cup R^+ =E_m\ \ \Longleftrightarrow\ \ R\succeq
K\ .
\end{equation}

\item[\rm(ii)]
For all $T',T''\in\mathbf{m}-\{B\}$, it holds
$(T')^+\cup(T'')^+\neq E_m$.

\item[\rm(iii)]
\begin{itemize}
\item[$\bullet$]

The subposet
\begin{equation}
\label{eq:34} \mathcal{O}(B):=\{T\in\mathcal{T}(\mathcal{L},B):\
B^+\cup T^+=E_m\}=\bigcap_{e\in B^-}\mathcal{T}_e^+
\end{equation}
is an order filter in the tope poset $\mathcal{T}(\mathcal{L},B)$,
with
\begin{equation*}
\bmin\mathcal{O}(B)=\mathcal{G}(B)\ ,
\end{equation*}
where the antichain $\mathcal{G}(B)$ is defined by
\begin{equation}
\label{eq:32} \mathcal{G}(B):=
\bigl\{T\in\mathcal{T}(\mathcal{L},B):\ T\in\bmax^+(\mathcal{T}),\
B^+\cup T^+=E_m\bigr\}\ .
\end{equation}

Furthermore, if $\mathcal{M}$ is totally cyclic, then it holds
\begin{equation*}
\mathcal{O}(B)=\conv_{\mathtt{T}}\bigl(\mathcal{G}(B)\bigr)\ .
\end{equation*}
\item[$\bullet$]
The union $\bigcup_{B\in\bmax^+(\mathcal{T})}\mathcal{O}(B)$
covers the set of topes $\mathcal{T}$ of\/ $\mathcal{M}$.
\item[$\bullet$]
For any topes $B',B''\in\bmax^+(\mathcal{T})$, we have
\begin{equation*}
|\mathcal{O}(B')\cap\mathcal{O}(B'')|>0\ \ \Longleftrightarrow\ \
|\mathcal{G}(B')\cap\mathcal{G}(B'')|>0\ .
\end{equation*}
\end{itemize}
\end{itemize}
\end{corollary}

\begin{proof}
(i) The uniqueness of the tope $K\in\bmax^+\bigl(\mathfrak{V}(\pmb{R})\bigr)$ such that $B^+$ $\cup K^+=E_m$,
and relation~(\ref{eq:30}) are discussed in the proof of Proposition~\ref{p:2} (substitute $T^0$, $T^l$ and
$T^m$ in that proof by $B$, $K$ and $-B$, respectively).

We have $|(-B)^+|=m-|B^+|$ and $|K^+|=|(-B)^+|+(m-\rho(K))$; hence
$|K^+|=2m-|B^+|-\rho(K)\leq c$, and~(\ref{eq:31}) follows.

(ii) This assertion is also inspired by the proof of Proposition~\ref{p:2}: if $\mathbf{m}$ $=(T^0,
T^1,\ldots, T^m)$, then for all $T',T''\in\mathbf{m}-\{B\}$, we have $(T')^+\cup(T'')^+\not\ni e$, where
$\{e\}:=\mathbf{S}(T^0,T^1)$.

(iii) The assertion follows from~(i).
\end{proof}

\begin{example}
Figure~{\rm\ref{fig:11}} depicts the Hasse diagram of a subposet
$\mathcal{O}(B)$ from Corollary~{\rm\ref{prop:36}(iii)} related to
a reorientation of an oriented matroid realized by the hyperplane
arrangement of Figure~{\rm\ref{fig:7}}.

\begin{figure}[ht]
\begin{picture}(12,2.5)(-1,0.75)
\put(2,1){\oval(2,0.5)} \put(5,1){\oval(2,0.5)}
\put(8,1){\oval(2,0.5)} \put(2,2){\oval(2,0.5)}
\put(5,2){\oval(2,0.5)} \put(8,2){\oval(2,0.5)}
\put(5,3){\oval(2,0.5)}

\put(2,1){\makebox(0,0)[c]{\scriptsize $+++-+-$}}
\put(5,1){\makebox(0,0)[c]{\scriptsize $-++++-$}}
\put(8,1){\makebox(0,0)[c]{\scriptsize $-+++-+$}}
\put(2,2){\makebox(0,0)[c]{\scriptsize $+++---$}}
\put(5,2){\makebox(0,0)[c]{\scriptsize $-++-+-$}}
\put(8,2){\makebox(0,0)[c]{\scriptsize $-+++--$}}
\put(5,3){\makebox(0,0)[c]{\scriptsize $-++---$}}

\put(3.9,3){\makebox(0,0)[r]{\small $-B$}}

\put(2,1.25){\line(0,1){0.5}} \put(5,1.25){\line(0,1){0.5}}
\put(8,1.25){\line(0,1){0.5}} \put(5,2.25){\line(0,1){0.5}}
\put(5,1.75){\line(-4,-1){2.1}} \put(8,1.75){\line(-4,-1){2.1}}
\put(5,2.75){\line(-4,-1){2.1}} \put(5,2.75){\line(4,-1){2.1}}
\end{picture}
\caption{The Hasse diagram of the order filter $\mathcal{O}(B)$ in
the tope poset $\mathcal{T}(\mathcal{L},B)$, defined
by~(\ref{eq:34}), for the reorientation
$\mathcal{M}:={}_{-[1,2]}\mathcal{N}^0$; the oriented matroid
$\mathcal{N}^0$ is realized by the hyperplane arrangement of
Figure~\ref{fig:7}, $B:=+--+++$.  The antichain
$\mathcal{G}(B):=\bmin\mathcal{O}(B)$ is
$\{+++-+-,-++++-,-+++-+\}\subset\bmax^+(\mathcal{T})$.}
\label{fig:11}
\end{figure}
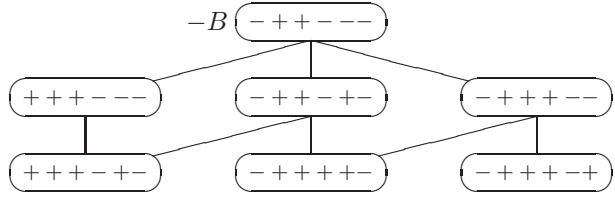
\end{example}

The following assertion (a proof of which we sketch in the Appendix) shows that Algorithm~{\rm\ref{alg:3}}
always constructs critical committees.

\begin{proposition}
\label{p:3} Let
$\mathcal{N}^0:=(E_m,\mathcal{L}^0)=(E_m,\mathcal{T}^0)$ be a
simple acyclic oriented matroid whose sets of covectors and of
topes are denoted by $\mathcal{L}^0$ and $\mathcal{T}^0$,
respectively.

Let $\mathbf{m}:=(R^0:=\mathrm{T}^{(+)}\ \precdot\ R^1\ \precdot\
\cdots\ \precdot\ R^m:=\mathrm{T}^{(-)})$ be a maximal chain in
the tope poset $\mathcal{T}(\mathcal{L}^0,\mathrm{T}^{(+)})$.

Let $s$ be an integer, $1\leq s\leq m$. Denote by $\mathcal{L}^s$
and $\mathcal{T}^s$ the sets of covectors and of topes,
respectively, of the reorientation
$\mathcal{N}^s:={}_{-[1,s]}\mathcal{N}^0$.

Let $\pmb{R}:=(T^0,T^1,\ldots,T^{2m-1},T^0)$ be a symmetric cycle
in the tope graph $\mathcal{T}^s(\mathcal{L}^s)$, such that
\begin{equation*}
T^k:={}_{-[1,s]}R^k\ ,\ \ \ 0\leq k\leq m\ .
\end{equation*}

Algorithm~{\rm\ref{alg:3}} builds the set
$\bmax^+\bigl(\mathfrak{V}(\pmb{R})\bigr)$ which is a critical
tope committee for $\mathcal{N}^s$.
\end{proposition}

\subsection{The Graph of Topes with Maximal Positive Parts}

Let $\mathcal{M}:=(E_m,\mathcal{L})=(E_m,\mathcal{T})$ be a simple oriented matroid. Choose in the graph
$\boldsymbol{\Gamma}(\mathcal{M})$, defined by~(\ref{eq:21}), the sub{\em graph of topes with
inclusion-maximal positive parts\/}
$\boldsymbol{\Gamma}_{\bmax}^+:=\boldsymbol{\Gamma}_{\bmax}^+(\mathcal{M})$ for which
\begin{equation}
\label{eq:33}
\begin{split}
\mathfrak{V}(\boldsymbol{\Gamma}_{\bmax}^+)\ \ &\ :=\ \
\bmax^+(\mathcal{T})\ ,\\
\{T^k,T^l\}\in\mathfrak{E}(\boldsymbol{\Gamma}_{\bmax}^+)\ \
&\Longleftrightarrow\ \ (T^k)^+\cup(T^l)^+=E_m\ .
\end{split}
\end{equation}

\begin{example}
The graph $\boldsymbol{\Gamma}_{\bmax}^+$ associated with a reorientation of an oriented matroid, that is
realized by the hyperplane arrangement of Figure~{\rm\ref{fig:7}}, is given in Figure~{\rm\ref{fig:12}}.

\begin{figure}[ht]
\begin{picture}(12,4.5)(-0.5,0.75)
\put(2,3){\oval(2,0.5)} \put(4,1){\oval(2,0.5)}
\put(4,5){\oval(2,0.5)} \put(7,1){\oval(2,0.5)}
\put(7,5){\oval(2,0.5)} \put(9,2){\oval(2,0.5)}
\put(9,4){\oval(2,0.5)}

\put(2,3){\makebox(0,0)[c]{\scriptsize $-+++-+$}}
\put(4,1){\makebox(0,0)[c]{\scriptsize $+++-+-$}}
\put(4,5){\makebox(0,0)[c]{\scriptsize $+--+++$}}
\put(7,1){\makebox(0,0)[c]{\scriptsize $--++++$}}
\put(7,5){\makebox(0,0)[c]{\scriptsize $++-+-+$}}
\put(9,2){\makebox(0,0)[c]{\scriptsize $-++++-$}}
\put(9,4){\makebox(0,0)[c]{\scriptsize $+-+-++$}}

\put(3,3){\line(1,-2){0.87}} \put(3,3){\line(1,2){0.87}}
\put(3,3){\line(5,1){5}} \put(4,1.25){\line(0,1){3.5}}
\put(8,2){\line(-4,3){3.65}} \put(4.2,1.25){\line(2,3){2.32}}
\put(5,1){\line(1,0){1}} \put(7,1.25){\line(0,1){3.5}}
\put(7.5,4.75){\line(1,-4){0.63}}
\put(7.75,4.75){\line(1,-1){0.5}} \put(9,2.25){\line(0,1){1.5}}
\end{picture}
\caption{The graph of topes with maximal positive parts
$\boldsymbol{\Gamma}_{\bmax}^+({}_{-[1,2]}\mathcal{N}^0)$, where
the oriented matroid $\mathcal{N}^0$ is realized by the hyperplane
arrangement of Figure~\ref{fig:7}.} \label{fig:12}
\end{figure}
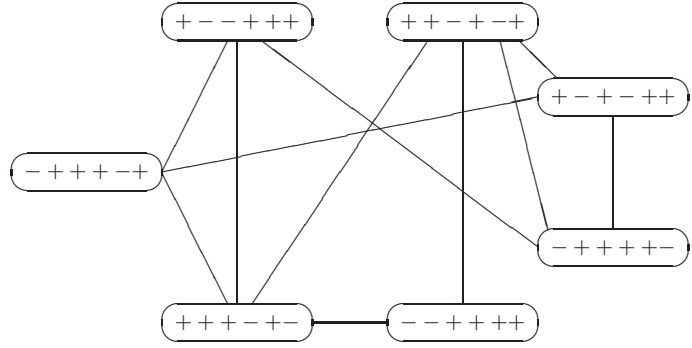
\end{example}

The graph $\boldsymbol{\Gamma}_{\bmax}^+$, defined
by~(\ref{eq:33}), is a direct generalization of the {\em graph of
maximal feasible subsystems\/} of an infeasible linear inequality
system which has been studied in
works~\cite{G-Thesis,G-PRIA,G-VINITI,G-1987,GNT}.

The {\em hypergraph of maximal feasible subsystems\/} of an infeasible linear inequality system is discussed,
e.g., in~\cite{Kh-PhDThesis,Kh-DMA,KhMR,MKh1}. An analogous construction for oriented matroids can be defined
in the following way: given a simple oriented matroid $\mathcal{M}=(E_m,\mathcal{T})$, the set of vertices of
the {\em hypergraph of topes with maximal positive parts
$\boldsymbol{\Xi}_{\bmax}^+:=\boldsymbol{\Xi}_{\bmax}^+(\mathcal{M})$} is the set $\bmax^+(\mathcal{T})$; a
subset $\mathcal{H}\subseteq\bmax^+(\mathcal{T})$ is a hyperedge of $\boldsymbol{\Xi}_{\bmax}^+$ iff
$\bigcup_{T\in\mathcal{H}}T^+=E_m$; thus the family of hyperedges $\mathfrak{E}(\boldsymbol{\Xi}_{\bmax}^+)$
of $\boldsymbol{\Xi}_{\bmax}^+$ is an order filter in the Boolean lattice of subsets of the set
$\bmax^+(\mathcal{T})$, with $\bmin\mathfrak{E}(\boldsymbol{\Xi}_{\bmax}^+)\supset
\mathfrak{E}(\boldsymbol{\Gamma}_{\bmax}^+)$.

A construction that is closely related to
$\boldsymbol{\Gamma}_{\bmax}^+(\mathcal{M})$ for realizable
coloopless simple oriented matroids $\mathcal{M}$ is the graph of
diagonals of a convex polytope. If $\pmb{P}$ is a convex polytope
with vertex set $\vvert\pmb{P}$ then a {\em diagonal\/} of
$\pmb{P}$ is a subset $D\subseteq\vvert\pmb{P}$ such that the
convex hull $\conv D$ of $D$ is not a proper face of $\pmb{P}$,
but $\conv(D-\{v\})$ lies in a proper face of $\pmb{P}$, for all
vertices $v\in D$, see, e.g., \cite{G-1985,G-Gu,Tam}; if,
furthermore, $\conv(V)$ is a face of $\pmb{P}$, for any proper
subset $V\subset D$, then $D$ is called an {\em empty simplex} (a
{\em missing face}) of $\pmb{P}$, see, e.g.,
\cite{K-Handbook,K-Aspects,N}. The {\em graph of diagonals\/} of
$\pmb{P}$ is defined as the Kneser graph of the family of
diagonals of $\pmb{P}$, see, e.g.~\cite[Definition~2.2.9]{M}.

Many properties of $\boldsymbol{\Gamma}_{\bmax}^+$, among which
the most important are connectedness and the existence of an odd
cycle, are inherited from the realizable case and lay the
foundation of graph-theoretic procedures of constructing  tope
committees of `high quality'.

\newpage

\subsubsection{General Properties of\/ $\boldsymbol{\Gamma}_{\bmax}^+$}

\begin{proposition}
Let $\mathcal{M}:=(E_m,\mathcal{L})=(E_m,\mathcal{T})$ be a simple
oriented matroid that is not acyclic.
\begin{itemize}
\item[\rm(i)]
The graph
$\boldsymbol{\Gamma}_{\bmax}^+:=\boldsymbol{\Gamma}_{\bmax}^+(\mathcal{M})$
is connected. The degree of every its vertex is at least two. Any
edge of $\boldsymbol{\Gamma}_{\bmax}^+$ is an edge of a cycle.
\item[\rm(ii)] If for any $2$-path $(R,B,S)$ in
$\boldsymbol{\Gamma}_{\bmax}^+$ there exist topes $R',S'\in\mathcal{T}$ $-\{-B\}$ such that
\begin{gather}
\nonumber (R')^+\subseteq R^+\ , \ \ \ (S')^+\subseteq S^+\ ,\\
\nonumber B^+\cup(R')^+=B^+\cup(S')^+=E_m \ ,\\ \label{eq:27}
|(R')^+\cap B^+\cap (S')^+|=0\ ,
\end{gather}
then $B$ is not a cutvertex in $\boldsymbol{\Gamma}_{\bmax}^+$.
\item[\rm(iii)] $\boldsymbol{\Gamma}_{\bmax}^+$
contains an odd cycle.
\end{itemize}
\end{proposition}

\begin{proof}
(i) Let $B$ and $R$ be any distinct topes from the set
$\bmax^+(\mathcal{T})$.

Let $\pmb{R}:=(T^0:=B,T^1,\ldots,T^k:=R,\ldots,T^{2m-1},T^0)$ be a symmetric cycle in the tope graph
$\mathcal{T}(\mathcal{L})$. By Proposition~\ref{p:2}, the set $\bmax^+\bigl(\mathfrak{V}(\pmb{R})\bigr)$ is
the set of vertices of an odd cycle $\pmb{C}$, defined in the following way: for
$T',T''\in\bmax^+\bigl(\mathfrak{V}(\pmb{R})\bigr)$, we have $\{T',T''\}\in\mathfrak{E}(\pmb{C})$ iff
$(T')^+\cup(T'')^+=E_m$.

Let $\phi:\
\mathfrak{V}(\pmb{C})\to\mathfrak{V}(\boldsymbol{\Gamma}_{\bmax}^+):=
\bmax^+(\mathcal{T})$ be any mapping such that
$T^+\subseteq\bigl(\phi(T)\bigr)^+$, for all
$T\in\mathfrak{V}(\pmb{C})$. This mapping is a graph homomorphism
from $\pmb{C}$ to $\boldsymbol{\Gamma}_{\bmax}^+$. Since $\pmb{C}$
is ($2$-)connected, there is a path in
$\boldsymbol{\Gamma}_{\bmax}^+$ between the vertices $\phi(T')$
and $\phi(T'')$, for all $T',T''\in\mathfrak{V}(\pmb{C})$. In
particular, there is a path in $\boldsymbol{\Gamma}_{\bmax}^+$
between $B$ and $R$ because $\phi(B)=B$ and $\phi(R)=R$.

Now suppose that $B^+\cup R^+=E_m$, that is,
$\{B,R\}\in\mathfrak{E}(\pmb{C})$. Let $\{B,T\}$ be the edge of
$\pmb{C}$ such that $T\neq R$; then
$\{B,R\},\{B,\phi(T)\}\in\mathfrak{E}(\boldsymbol{\Gamma}_{\bmax}^+)$,
where $\phi(T)\neq B$ and $\phi(T)\neq R$, therefore the degree of
$B$ in $\boldsymbol{\Gamma}_{\bmax}^+$ is greater than one.

Let $\pmb{D}$ denote the path in the cycle $\pmb{C}$ between the
vertices $B$ and $R$, such that $T\in\mathfrak{V}(\pmb{D})$. The
image of $\pmb{D}$ under the homomorphism $\phi$ is a connected
subgraph of $\boldsymbol{\Gamma}_{\bmax}^+$ whose set of edges
does not contain the edge $\{B,R\}=\{\phi(B),\phi(R)\}$. Hence,
the edge $\{B,R\}\in\mathfrak{E}(\boldsymbol{\Gamma}_{\bmax}^+)$
is an edge of a cycle.

(ii) Assume that $R$ and $S$ belong to different blocks of the
graph $\boldsymbol{\Gamma}_{\bmax}^+$, that is, $B$ is a
cutvertex.

Since, by the hypothesis of the assertion, it holds $B^+\cup (R')^+=B^+$ $\cup (S')^+=E_m$,
condition~(\ref{eq:27}) implies
\begin{equation*}
\mathbf{S}(B,-R')\subset\mathbf{S}(B,S')\ ,\ \ \
\mathbf{S}(B,-S')\subset\mathbf{S}(B,R')
\end{equation*}
and, as a consequence,
\begin{equation*}
S'\succ -R'\ ,\ \ \ R'\succ -S'
\end{equation*}
in the tope poset $\mathcal{T}(\mathcal{L},B)$. This implies that in the tope graph
$\mathcal{T}(\mathcal{L})$ there exists a symmetric cycle $\pmb{R}$ such that
$\{R',B,S'\}\subset\mathfrak{V}(\pmb{R})$. Let $R''$ and $S''$ be the topes with $(R')^+\subseteq(R'')^+$ and
$(S')^+\subseteq(S'')^+$, such that $R'',S''$ $\in\bmax^+\bigl(\mathfrak{V}(\pmb{R})\bigr)$.

Let $\pmb{C}$ and $\phi:\
\mathfrak{V}(\pmb{C})\to\mathfrak{V}(\boldsymbol{\Gamma}_{\bmax}^+)$
be an odd cycle and a graph homomorphism, respectively, which were
defined in the proof of assertion~(i), with $\phi(R''):=R$ and
$\phi(S''):=S$.

By Proposition~\ref{p:2}, the sets $\{R'',B\}$ and $\{B,S''\}$ are
edges of $\pmb{C}$ and, as a consequence,
$\bigl(\phi(R''),\phi(B),$ $\phi(S'')\bigr)=(R,B,S)$ is a $2$-path
in $\boldsymbol{\Gamma}_{\bmax}^+$. Let $\pmb{D}$ denote the path
in the cycle $\pmb{C}$ between the vertices $R''$ and $S''$, that
does not contain $B$. The image of $\pmb{D}$ under the
homomorphism $\phi$ is a connected subgraph of
$\boldsymbol{\Gamma}_{\bmax}^+$ whose set of vertices does not
contain the tope $B=\phi(B)$. Hence, there is a cycle in
$\boldsymbol{\Gamma}_{\bmax}^+$ such that $\{R,B,S\}$ is a subset
of its vertices. This contradicts our assumption that $B$ is a
cutvertex in $\boldsymbol{\Gamma}_{\bmax}^+$.

(iii) Let $\pmb{R}:=(T^0,T^1,\ldots,T^{2m-1},T^0)$ be a symmetric
cycle in the tope graph $\mathcal{T}(\mathcal{L})$, such that
$T^0\in\bmax^+(\mathcal{T})$. Let $\pmb{G}$ be the graph defined
by~(\ref{eq:8}). Recall that the set
$\bmax^+\bigl(\mathfrak{V}(\pmb{R})\bigr)$ is the set of vertices
$\mathfrak{V}(\pmb{C})$ of an odd cycle $\pmb{C}$ in $\pmb{G}$,
see Proposition~\ref{p:2}. The set $\mathfrak{E}(\pmb{C})$ of
edges of $\pmb{C}$ is described by~(\ref{eq:10}).

Assume that $\boldsymbol{\Gamma}_{\bmax}^+$ is bipartite, with
partition classes $\mathfrak{V}'$ and $\mathfrak{V}''$. Suppose
that $T^0\in\mathfrak{V}'$. Then, for a homomorphism
\begin{equation*}
\phi:\
\mathfrak{V}(\pmb{C}):=\bmax^+\bigl(\mathfrak{V}(\pmb{R})\bigr)\ \
\to\ \
\mathfrak{V}(\boldsymbol{\Gamma}_{\bmax}^+):=\bmax^+(\mathcal{T})\
\end{equation*}
from $\pmb{C}$ to $\boldsymbol{\Gamma}_{\bmax}^+$, such that
$T^+\subseteq\bigl(\phi(T)\bigr)^+$ for all
$T\in\mathfrak{V}(\pmb{C})$, we have
\begin{gather}
\nonumber \phi(T^0)=T^0\in\mathfrak{V}'\ ,\ \ \{T^0,T^{k_{d+1}}\},
\{T^0,T^{k_d}\}\in\mathfrak{E}(\pmb{C})\ \ \ \Longrightarrow\ \ \
\phi(T^{k_{d+1}})\in\mathfrak{V}''\ ,\\ \label{eq:28}
\phi(T^{k_d})\in\mathfrak{V}''\ ;\\ \nonumber
\phi(T^{k_{d+1}})\in\mathfrak{V}''\ ,\ \
\{T^{k_1},T^{k_{d+1}}\}\in\mathfrak{E}(\pmb{C})\ \ \
\Longrightarrow\ \ \ \phi(T^{k_1})\in\mathfrak{V}'\ ;\\ \nonumber
\phi(T^{k_1})\in\mathfrak{V}'\ ,\ \
\{T^{k_1},T^{k_{d+2}}\}\in\mathfrak{E}(\pmb{C})\ \ \
\Longrightarrow\ \ \ \phi(T^{k_{d+2}})\in\mathfrak{V}''\ ;\\
\nonumber \vdots\\ \nonumber \phi(T^{k_{2d}})\in\mathfrak{V}''\ ,\
\ \{T^{k_d},T^{k_{2d}}\}\in\mathfrak{E}(\pmb{C})\ \ \
\Longrightarrow\\ \label{eq:29} \phi(T^{k_d})\in\mathfrak{V}'.
\end{gather}
Since~(\ref{eq:29}) contradicts~(\ref{eq:28}),
$\boldsymbol{\Gamma}_{\bmax}^+$ is not bipartite; as a
consequence, it contains an odd cycle, see, e.g.,
~\cite[Proposition~1.6.1]{Diestel}.
\end{proof}

Suppose that $\boldsymbol{\Gamma}_{\bmax}^+$ has no cutvertices,
that is, $\boldsymbol{\Gamma}_{\bmax}^+$ is $2$-connected. Since
$\boldsymbol{\Gamma}_{\bmax}^+$ contains an odd cycle then the
general properties of $2$-connected graphs~\cite[\S{}5.4]{Ore}
imply that every vertex of $\boldsymbol{\Gamma}_{\bmax}^+$ is
contained in an odd cycle.

\subsubsection{The Neighborhood Complex of\/ $\boldsymbol{\Gamma}_{\bmax}^+$}
Let $\mathcal{M}:=(E_m,\mathcal{L})=(E_m,\mathcal{T})$ be a simple
oriented matroid that is not acyclic. Recall that the neighborhood
of a vertex $B$ in
$\boldsymbol{\Gamma}_{\bmax}^+:=\boldsymbol{\Gamma}_{\bmax}^+(\mathcal{M})$
is the set $\{T\in\bmax^+(\mathcal{T}):\ B^+\cup T^+=E_m\}$.
Equivalently, the neighborhood of $B$ is the antichain
$\mathcal{G}(B)$ in the tope poset $\mathcal{T}(\mathcal{L},B)$,
defined by~(\ref{eq:32}).

According to the {\em Folkman-Lawrence Topological Representation
Theorem}~\cite[\S{}5.2]{BLSWZ}, one can consider a representation
of $\mathcal{M}$ by an arrangement of oriented pseudospheres
\begin{equation*}
\{S_e:\ e\in E_m\}
\end{equation*}
lying on the standard $(r(\mathcal{M})-1)$-dimensional sphere
$\mathbb{S}^{r(\mathcal{M})-1}$, where $r(\mathcal{M})$ denotes
the rank of $\mathcal{M}$, see~\cite{ERW}.

Recall that the simplicial complex of acyclic subsets of $E_m$,
denoted by $\Delta_{\acyclic}(\mathcal{M})$, is the nerve of the
family
\begin{equation}
\label{eq:23} \{S_e^+:\ e\in E_m\}\ ,
\end{equation}
where $S_e^+$ denotes the open positive hemisphere corresponding
to the pseudosphere $S_e$, that is, the positive side of $S_e$,
see~\cite{ERW}.

$\bullet$\ Suppose $\mathcal{M}$ is totally cyclic. Recall that in
this case the complex $\Delta_{\acyclic}$ is homotopy equivalent
to $\mathbb{S}^{r(\mathcal{M})-1}$ because the union of the sets
from~(\ref{eq:23}) is an open cover of the sphere by subspaces
whose nonempty intersections are contractible:
\begin{equation*}
\bigcup_{e\in E_m}S_e^+=\mathbb{S}^{r(\mathcal{M})-1}\ ,
\end{equation*}
see~\cite{ERW}. Combinatorial homotopy is discussed, e.g.,
in~\cite[\S{}10]{B-TM}.

By Corollary~\ref{prop:36}(iii), the neighborhood complex
$\NC(\boldsymbol{\Gamma}_{\bmax}^+)$ of
$\boldsymbol{\Gamma}_{\bmax}^+$ is the nerve of the family of open
subspaces
\begin{equation}
\label{eq:22} \Bigl\{\bigcap_{e\in B^-}S_e^+:\
B\in\bmax^+(\mathcal{T})\Bigr\}\
\end{equation}
such that every their nonempty intersection is contractible.

The {\em Nerve Theorem}~\cite[Theorem~10.7]{B-TM} implies that
$\NC(\boldsymbol{\Gamma}_{\bmax}^+)$ is homotopy equivalent to the
subspace of $\mathbb{S}^{r(\mathcal{M})-1}$ covered by
family~(\ref{eq:22}).

$\bullet$\ Suppose $\mathcal{M}$ is neither acyclic nor totally cyclic. It is shown in~\cite{ERW} that there
exists a unique non-negative covector $F\in\mathcal{L}$ with inclusion-maximum positive part. Denote by
$\boldsymbol{\Gamma}_{\bmax}^+(\mathcal{M}\backslash F^+)$ the graph of topes with maximal positive parts,
which is associated with the (totally cyclic) deletion $\mathcal{M}\backslash F^+$. The graph
$\boldsymbol{\Gamma}_{\bmax}^+(\mathcal{M})$ is isomorphic to
$\boldsymbol{\Gamma}_{\bmax}^+(\mathcal{M}\backslash F^+)$.

\section{Committees, Halfspaces, and Relatively Blocking
Elements} \label{section:5} Let $\mathcal{M}$ {$=(E_m,\mathcal{T})$ be a simple oriented matroid with set of
topes $\mathcal{T}$.

Denote by $\pmb{I}_{\frac{1}{2}}(\mathcal{T}_1^+,\ldots,\mathcal{T}_m^+)$ the family of all set-theoretic
committees for the family of positive halfspaces $\{\mathcal{T}_1^+,\ldots,\mathcal{T}_m^+\}$: by definition,
a set $\mathcal{P}\subset\mathcal{T}$ is a {\em committee\/} for the family of positive halfspaces
$\{\mathcal{T}_1^+,\ldots,\mathcal{T}_m^+\}$ iff it holds
\begin{equation*}
|\mathcal{P}\cap\mathcal{T}_e^+|>\frac{|\mathcal{P}|}{2}\ ,
\end{equation*}
for all $e\in E_m$. For every $e\in E_m$, we have $|\mathcal{T}_e^+|=\tfrac{|\mathcal{T}|}{2}$; the way of
computing the cardinality of the tope set $\mathcal{T}$ is
well-known~\cite[Theorem~4.6.1]{BLSWZ},\cite{LasV1975a,ThZ1975b,ThZ1975a}.

We have
\begin{equation*} \mathbf{K}^{\ast}(\mathcal{M})=
\pmb{I}_{\frac{1}{2}}(\mathcal{T}_1^+,\ldots,\mathcal{T}_m^+)\ ,
\end{equation*}
that is, a set $\mathcal{K}^{\ast}\subset\mathcal{T}$ is a tope committee for $\mathcal{M}$ iff
$\mathcal{K}^{\ast}$ is a set-theoretic committee for the family
$\{\mathcal{T}_1^+,\ldots,\mathcal{T}_m^+\}$.

Let $\mathbb{B}(\mathcal{T})$ denote the Boolean lattice of all subsets of $\mathcal{T}$. For any $e\in E_m$,
we relate to the set $\mathcal{T}_e^+$ the element
$\upsilon_e:=\bigvee_{T\in\mathcal{T}_e^+}T\in\mathbb{B}(\mathcal{T})$, the join of those atoms of
$\mathbb{B}(\mathcal{T})$ that compose $\mathcal{T}_e^+$.

Consider every tope committee for $\mathcal{M}$ as an element of the Boolean lattice
$\mathbb{B}(\mathcal{T})$. Then $\mathbf{K}^{\ast}(\mathcal{M})$ is precisely the subposet
$\mathbf{I}_{\frac{1}{2}}\bigl(\mathbb{B}(\mathcal{T}), \{\upsilon_1,\ldots,\upsilon_m\}\bigr)$ of all {\em
relatively \mbox{$\tfrac{1}{2}$-blocking elements}\/} for the antichain $\{\upsilon_1,\ldots,\upsilon_m\}$ in
$\mathbb{B}(\mathcal{T})$, with respect to the poset rank function on $\mathbb{B}(\mathcal{T})$. Relative
blocking in posets is discussed in~\cite{AM}. The antichain
$\mathfrak{y}_{\tfrac{1}{2}}(\mathbb{B}(\mathcal{T}),
\{\upsilon_1,\ldots,\upsilon_m\}):=\bmin\mathbf{I}_{\frac{1}{2}}\bigl(\mathbb{B}(\mathcal{T}),$
$\{\upsilon_1,\ldots,\upsilon_m\}\bigr)$, called in~\cite{AM} the {\em relative $\tfrac{1}{2}$-blocker of\/
$\{\upsilon_1,\ldots,$ $\upsilon_m\}$ in $\mathbb{B}(\mathcal{T})$}, is the family of all minimal committees
for $\mathcal{M}$.

For $k\in[1,|\mathcal{T}|-1]$, the subposet $\mathbf{I}_{\frac{1}{2},k}\bigl(\mathbb{B}(\mathcal{T}),
\{\upsilon_1,\ldots,\upsilon_m\}\bigr)$ of elements of rank $k$ from
$\mathbf{I}_{\frac{1}{2}}\bigl(\mathbb{B}(\mathcal{T}), \{\upsilon_1,\ldots,\upsilon_m\}\bigr)$ is
\begin{equation*}
\mathbf{I}_{\frac{1}{2},k}\bigl(\mathbb{B}(\mathcal{T}),
\{\upsilon_1,\ldots,\upsilon_m\}\bigr)=\mathbb{B}(\mathcal{T})^{(k)}\cap \bigcap_{e\in
E_m}\mathfrak{F}\!\left(\mathfrak{I}(\upsilon_e)\ \cap\
\mathbb{B}(\mathcal{T})^{(\lceil(k+1)/2\rceil)}\right)\ ,
\end{equation*}
according to~\cite[Proposition~5.1(ii)]{AM}. Layers $\mathbf{I}_{\frac{1}{2},k}\bigl(\mathbb{B}(\mathcal{T}),
\{\upsilon_1,\ldots,\upsilon_m\}\bigr)$ are discussed in~\cite{AM-P}, and Farey subsequences appearing in
their analysis are considered in~\cite{AM-Int,AM-P,AM}.

\section*{Appendix}
We sketch here a proof of Proposition~\ref{p:3}.

\begin{sketch}
Argue by induction on $i$, $1\leq i\leq s$.

The following observation is useful:

{\sl Claim.} {\em For any $t'$ and $t''$ such that $0\leq t'<t''\leq m$, it holds
\begin{equation*}
\begin{split}
\bigl(({}_{-[1,i]}R^{t'})\!\mid_{[1,i]}\bigr)^+\ &\subseteq\
\bigl(({}_{-[1,i]}R^{t''})\!\mid_{[1,i]}\bigr)^+\ ,\\ \bigl(({}_{-[1,i]}R^{t'})\!\mid_{[i+1,m]}\bigr)^+\
&\supseteq\ \bigl(({}_{-[1,i]}R^{t''})\!\mid_{[i+1,m]}\bigr)^+\ .
\end{split}
\end{equation*}
For any $t'$ and $t''$ such that $1\leq t'<t''\leq m$, it holds
\begin{equation*}
\begin{split}
\bigl(({}_{-[i,m]}R^{t'})\!\mid_{[1,i]}\bigr)^+\ &\supseteq\
\bigl(({}_{-[i,m]}R^{t''})\!\mid_{[1,i]}\bigr)^+\ ,\\ \bigl(({}_{-[i,m]}R^{t'})\!\mid_{[i+1,m]}\bigr)^+\
&\subseteq\ \bigl(({}_{-[i,m]}R^{t''})\!\mid_{[i+1,m]}\bigr)^+\ .
\end{split}
\end{equation*}
} \hfill $\vartriangleright\vartriangleleft$

Assign to every tope $R^i\in\{R^1,\ldots,R^m\}$ the label $\mathfrak{l}_i$ defined by~(\ref{eq:14}).

$\quad$

Suppose $i=1$.

\begin{itemize}
\item[$\bullet$]
If $\mathfrak{l}_1=i$ then we have $\mathcal{K}_i{}^{\ast}=\{{}_{-[1,i]}R^1\}=\{\mathrm{T}^{(+)}\}$ because
at Step~{\it 03\/} the algorithm builds the set
$\mathcal{K}_i{}^{\ast}=\{{}_{-[1,i]}R^1,{}_{-[1,i]}R^0,{}_{-[i,m]}R^1\}$, and the pair of opposites
$\{{}_{-[1,i]}R^0,{}_{-[i,m]}R^1\}$ is removed at Steps~{\it 04-05}.
\item[$\bullet$]
If $\mathfrak{l}_m=i$ then we have $\mathcal{K}_i{}^{\ast}=\{{}_{-[i,m]}R^m\}=\{\mathrm{T}^{(+)}\}$ because
at Step~{\it 03\/} the algorithm builds the set
$\mathcal{K}_i{}^{\ast}=\{{}_{-[i,m]}R^m,{}_{-[1,i]}R^0,{}_{-[1,i]}R^m\}$, and the pair of opposites
$\{{}_{-[1,i]}R^0,{}_{-[1,i]}R^m\}$ is removed at Steps~{\it 04-05}.
\item[$\bullet$]
If $\mathfrak{l}_j=i$, for some $j$, $1<j<m$, then we have
$\mathcal{K}_i{}^{\ast}=\{{}_{-[1,i]}R^0,{}_{-[1,i]}R^k,$ ${}_{-[i,m]}R^k\}$. Since
\begin{gather*}
({}_{-\mathfrak{l}_m}(-R^m))(i)=R^0(i)=R^1(i)=+\ ,\ \ \
({}_{-\mathfrak{l}_m}(-R^m))^+\subsetneqq(R^0)^+\supsetneqq(R^1)^+\ ,\\ R^{j-1}(i)=+\ ,\ \ \
R^j(i)=R^{j+1}(i)=-\ ,\ \ \ (R^{j-1})^+\supsetneqq (R^j)^+\supsetneqq(R^{j+1})^+\ ,
\end{gather*}
and
\begin{multline*}
({}_{-\mathfrak{l}_{j-1}}(-R^{j-1}))(i)=({}_{-\mathfrak{l}_j}(-R^j))(i)=-\ ,\ \ \
({}_{-\mathfrak{l}_{j+1}}(-R^{j+1}))(i)=+\ ,\\
({}_{-\mathfrak{l}_{j-1}}(-R^{j-1}))^+\subsetneqq({}_{-\mathfrak{l}_j}(-R^j))^+
\subsetneqq({}_{-\mathfrak{l}_{j+1}}(-R^{j+1}))^+\ ,
\end{multline*}
we have
\begin{gather*}
({}_{-[i,m]}R^m)^+\subsetneqq(R^0)^+\supsetneqq(R^1)^+\ ,
\\ ({}_{-[1,i]}R^{j-1})^+\subsetneqq({}_{-[1,i]}R^j)^+\supsetneqq
({}_{-[1,i]}R^{j+1})^+\\ ({}_{-[i,m]}R^{j-1})^+\subsetneqq({}_{-[i,m]}R^j)^+
\supsetneqq({}_{-[i,m]}R^{j+1})^+ ;
\end{gather*}
therefore
\begin{equation*}
(R^0)^+,({}_{-[1,i]}R^j)^+,({}_{-[i,m]}R^j)^+\in\bmax\{({}_{-[1,i]}R^k)^+:\ 0\leq k\leq 2m-1\}\ ,
\end{equation*}
see Remark~\ref{p:8}(ii).

One can show by means of a similar argument that for every $R\in\mathfrak{V}(\pmb{R})$ such that
$R\not\in\{R^0,R^j,{}_{-\mathfrak{l}_j}(-R^j)\}$, it holds $R^+\not\in\bmax\{({}_{-[1,i]}R^k)^+:\ 0\leq k\leq
2m-1\}$.
\end{itemize}

Thus,
\begin{equation*}
\{K^+:\ K\in\mathcal{K}_i{}^{\ast}\}=\bmax\{({}_{-[1,i]}R^k)^+:\ 0\leq k\leq 2m-1\}\ .
\end{equation*}

If $i=s$ then we are done.

Suppose $i>1$.

By the induction hypothesis, we have
\begin{equation*}
\{K^+:\ K\in\mathcal{K}_{i-1}{}^{\ast}\}=\bmax\{({}_{-[1,i-1]}R^k)^+:\ 0\leq k\leq 2m-1\}\ .
\end{equation*}

Suppose that $\mathfrak{l}_j=i$, for some $j$, $1\leq j\leq m$.

Consider the set
\begin{equation*}
\overline{\mathcal{K}_i{}^{\ast}}:=\{{}_{-i}K:\ K\in\mathcal{K}_{i-1}{}^{\ast}\}\ \dot\cup\
\{{}_{-[1,i]}R^j,{}_{-[i,m]}R^j\}\ .
\end{equation*}

\begin{itemize}
\item[$\bullet$]
Suppose $j=1$.

$\diamond$ Suppose that
\begin{equation}
\label{eq:15} {}_{-[1,i]}R^0\not\in\overline{\mathcal{K}_i{}^{\ast}}\ ,\ \ \
{}_{-[i,m]}R^{j+1}\not\in\overline{\mathcal{K}_i{}^{\ast}}\ .
\end{equation}

For $i>1$, we only describe the induction step for the case where $j=1$ and~(\ref{eq:15}) holds. Analysis of
other situations is completely analogous.

Condition~(\ref{eq:15}) implies that
\begin{equation*}
{}_{-[1,i]}R^{m-1}\not\in\overline{\mathcal{K}_i{}^{\ast}}\ ,\ \ \
{}_{-[1,i]}R^m\not\in\overline{\mathcal{K}_i{}^{\ast}}\ ,\ \ \
{}_{-[1,i]}R^{j+1}\not\in\overline{\mathcal{K}_i{}^{\ast}}\ ,
\end{equation*}
and it holds
\begin{equation*}
({}_{-[i-1,m]}R^m)^+\supsetneqq ({}_{-[1,i-1]}R^0)^+\supsetneqq({}_{-[1,i-1]}R^j)^+
\supsetneqq({}_{-[1,i-1]}R^{j+1})^+\ ;
\end{equation*}
\begin{multline*}
({}_{-[1,i-1]}R^{m-1})^+ \subsetneqq({}_{-[1,i-1]}R^m)^+= ({}_{-[i-1,m]}R^j)^+\\
\subsetneqq({}_{-[i-1,m]}R^{j+1})^+ \subsetneqq({}_{-[i-1,m]}R^{j+2})^+\ .
\end{multline*}
Note that
\begin{gather*}
({}_{-[1,i-1]}R^0)(e)=-\ ,\ 1\leq e\leq i-1\ ;\ \ \ ({}_{-[1,i-1]}R^0)(e)=+\ ,\  i\leq e\leq m\ ;\\
({}_{-[1,i-1]}R^j)(e)=-\ ,\ 1\leq e\leq i\ ;\ \ \ ({}_{-[1,i-1]}R^j)(e)=+\ ,\  i+1\leq e\leq m\ ;\\
({}_{-[1,i-1]}R^{j+1})(e)=-\ ,\ 1\leq e\leq i\ .
\end{gather*}
As a consequence, we have
\begin{equation*}
({}_{-[1,i]}R^0)^+\subsetneqq({}_{-[1,i]}R^j)^+ \supsetneqq({}_{-[1,i]}R^{j+1})^+\ ,
\end{equation*}
that is,
\begin{equation}
\label{eq:16} ({}_{-[1,i]}R^j)^+\in\bmax\{({}_{-[1,i]}R^k)^+:\ 0\leq k\leq 2m-1\}\ ,
\end{equation}
see Remark~\ref{p:8}(ii).

Since
\begin{gather*}
({}_{-[1,i-1]}R^{m-1})(e)=-\ ,\  i\leq e\leq m\ ;\\ ({}_{-[1,i-1]}R^m)(e)=+\ ,\ 1\leq e\leq i-1\ ;\ \ \
({}_{-[1,i-1]}R^m)(e)=-\ ,\  i\leq e\leq m\ ,
\end{gather*}
and $({}_{-[i,m]}R^{j+1})^+=(-{}_{-[1,i]}R^j)^+\in\bmin\{({}_{-[1,i]}R^k)^+:\ 0\leq k\leq 2m-1\}$,
by~(\ref{eq:16}), we obtain
\begin{equation*}
({}_{-[1,i]}R^{m-1})^+ \subsetneqq({}_{-[1,i]}R^m)^+= ({}_{-[i,m]}R^j)^+\\ \supsetneqq({}_{-[i,m]}R^{j+1})^+\
,
\end{equation*}
that is,
\begin{equation*}
({}_{-[i,m]}R^j)^+\in\bmax\{({}_{-[1,i]}R^k)^+:\ 0\leq k\leq 2m-1\}\ ,
\end{equation*}
by Remark~\ref{p:8}(ii).

Note that for all $k$, $1\leq k\leq m$, it holds $({}_{-[1,i]}R^k)(i)=+$, cf.~Remark~\ref{p:8}(i).

Let ${}_{-i}K\in\overline{\mathcal{K}_i{}^{\ast}}- \{{}_{-[1,i]}R^j,{}_{-[i,m]}R^j\}$, that is,
$K\in\mathcal{K}_{i-1}{}^{\ast}$. Let $(R',K,R'')$ be a $2$-path in the cycle
$({}_{-[1,i-1]}R^0,{}_{-[1,i-1]}R^1,\ldots, {}_{-[1,i-1]}R^{2m-1},$ ${}_{-[1,i-1]}R^0)$; by
Remark~\ref{p:8}(ii), the vertices of the path satisfy $(R')^+\subsetneqq K^+ \supsetneqq (R'')^+$.
Remark~\ref{p:8}(i) implies that $R'(i)=K(i)=R''(i)$. As a consequence, the equality
$({}_{-i}R')(i)=({}_{-i}K)(i)=({}_{-i}R'')(i)$ holds as well, and we have
$({}_{-i}R')^+\subsetneqq({}_{-i}K)^+\supsetneqq({}_{-i}R'')^+$. Thus, we have $\{({}_{-i}K)^+:\
K\in\mathcal{K}_{i-1}{}^{\ast}\}\subset \bmax\{({}_{-[1,i]}R^k)^+:\ 0\leq k\leq 2m-1\}$.

In a similar manner, one can show that for any $R\in\{{}_{-[1,i]}R^0,{}_{-[1,i]}R^1,$
$\ldots,{}_{-[1,i]}R^{2m-1}\}$ such that ${}_{-[1,i]}R\not\in\overline{\mathcal{K}_i{}^{\ast}}$ it holds
$({}_{-[1,i]}R)^+\not\in\bmax\{({}_{-[1,i]}R^k)^+:\ 0\leq k\leq 2m-1\}$.

The algorithm builds the set
\begin{equation*}
\mathcal{K}_i{}^{\ast}=\overline{\mathcal{K}_i{}^{\ast}}\ ;
\end{equation*}
we have seen that
\begin{equation*}
\{K^+:\ K\in\mathcal{K}_i{}^{\ast}\}=\bmax\{({}_{-[1,i]}R^k)^+:\ 0\leq k\leq 2m-1\}\ .
\end{equation*}

$\diamond$ If
\begin{equation*}
{}_{-[1,i]}R^0\in\overline{\mathcal{K}_i{}^{\ast}}\ ,\ \ \
{}_{-[i,m]}R^{j+1}\not\in\overline{\mathcal{K}_i{}^{\ast}}\ ,
\end{equation*}
then
\begin{equation*}
\begin{split}
\mathcal{K}_i{}^{\ast}&=\overline{\mathcal{K}_i{}^{\ast}}\ - \
\{{}_{-[1,i]}R^0,{}_{-[i,m]}R^j\}\\&=\bigl(\{{}_{-i}K:\
K\in\mathcal{K}_{i-1}{}^{\ast}\}-\{{}_{-[1,i]}R^0\}\bigr)\ \dot\cup\ \{{}_{-[1,i]}R^j\}\ .
\end{split}
\end{equation*}

$\diamond$ If
\begin{equation*}
{}_{-[1,i]}R^0\in\overline{\mathcal{K}_i{}^{\ast}}\ ,\ \ \
{}_{-[i,m]}R^{j+1}\in\overline{\mathcal{K}_i{}^{\ast}}\ ,
\end{equation*}
then
\begin{equation*}
\begin{split}
\mathcal{K}_i{}^{\ast}&=\overline{\mathcal{K}_i{}^{\ast}}\ - \
\{{}_{-[1,i]}R^0,{}_{-[1,i]}R^j,{}_{-[i,m]}R^j,{}_{-[i,m]}R^{j+1}\}\\&= \{{}_{-i}K:\
K\in\mathcal{K}_{i-1}{}^{\ast}\}\ -\ \{{}_{-[1,i]}R^0,{}_{-[i,m]}R^{j+1}\} \ .
\end{split}
\end{equation*}

If $i=s$ then we are done.

\item[$\bullet$] Suppose $j=m$.

$\diamond$ If
\begin{equation*}
{}_{-[1,i]}R^0\not\in\overline{\mathcal{K}_i{}^{\ast}}\ ,\ \ \
{}_{-[1,i]}R^{j-1}\not\in\overline{\mathcal{K}_i{}^{\ast}}\ ,
\end{equation*}
then
\begin{equation*}
\mathcal{K}_i{}^{\ast}=\overline{\mathcal{K}_i{}^{\ast}}\ .
\end{equation*}

$\diamond$ If
\begin{equation*}
{}_{-[1,i]}R^0\in\overline{\mathcal{K}_i{}^{\ast}}\ ,\ \ \
{}_{-[1,i]}R^{j-1}\not\in\overline{\mathcal{K}_i{}^{\ast}}\ ,
\end{equation*}
then
\begin{equation*}
\begin{split}
\mathcal{K}_i{}^{\ast}&=\overline{\mathcal{K}_i{}^{\ast}}\ - \ \{{}_{-[1,i]}R^0,{}_{-[1,i]}R^j\}\\&=
\bigl(\{{}_{-i}K:\ K\in\mathcal{K}_{i-1}{}^{\ast}\}-\{{}_{-[1,i]}R^0\}\bigr)\ \dot\cup\ \{{}_{-[i,m]}R^j\} \
.
\end{split}
\end{equation*}

$\diamond$ If
\begin{equation*}
{}_{-[1,i]}R^0\in\overline{\mathcal{K}_i{}^{\ast}}\ ,\ \ \
{}_{-[1,i]}R^{j-1}\in\overline{\mathcal{K}_i{}^{\ast}}\ ,
\end{equation*}
then
\begin{equation*}
\begin{split}
\mathcal{K}_i{}^{\ast}&=\overline{\mathcal{K}_i{}^{\ast}}\ - \
\{{}_{-[1,i]}R^0,{}_{-[1,i]}R^{j-1},{}_{-[1,i]}R^m,{}_{-[i,m]}R^j\}\\&= \{{}_{-i}K:\
K\in\mathcal{K}_{i-1}{}^{\ast}\}\ - \ \{{}_{-[1,i]}R^0,{}_{-[1,i]}R^{j-1}\} \ .
\end{split}
\end{equation*}

\item[$\bullet$] Suppose $1<j<m$.

$\diamond$ If
\begin{equation*}
{}_{-[1,i]}R^{j-1}\not\in\overline{\mathcal{K}_i{}^{\ast}}\ ,\ \ \
{}_{-[i,m]}R^{j+1}\not\in\overline{\mathcal{K}_i{}^{\ast}}\ ,
\end{equation*}
then
\begin{equation*}
\mathcal{K}_i{}^{\ast}=\overline{\mathcal{K}_i{}^{\ast}}\ .
\end{equation*}

$\diamond$ If
\begin{equation*}
{}_{-[1,i]}R^{j-1}\in\overline{\mathcal{K}_i{}^{\ast}}\ ,\ \ \
{}_{-[i,m]}R^{j+1}\not\in\overline{\mathcal{K}_i{}^{\ast}}\ ,
\end{equation*}
then
\begin{equation*}
\begin{split}
\mathcal{K}_i{}^{\ast}&=\overline{\mathcal{K}_i{}^{\ast}}\ - \
\{{}_{-[1,i]}R^{j-1},{}_{-[i,m]}R^j\}\\&=\bigl(\{{}_{-i}K:\
K\in\mathcal{K}_{i-1}{}^{\ast}\}-\{{}_{-[1,i]}R^{j-1}\}\bigr)\ \dot\cup\ \{{}_{-[1,i]}R^j\} \ .
\end{split}
\end{equation*}

$\diamond$ If
\begin{equation*}
{}_{-[1,i]}R^{j-1}\in\overline{\mathcal{K}_i{}^{\ast}}\ ,\ \ \
{}_{-[i,m]}R^{j+1}\in\overline{\mathcal{K}_i{}^{\ast}}\ ,
\end{equation*}
then
\begin{equation*}
\begin{split}
\mathcal{K}_i{}^{\ast}&=\overline{\mathcal{K}_i{}^{\ast}}\ - \
\{{}_{-[1,i]}R^{j-1},{}_{-[1,i]}R^j,{}_{-[i,m]}R^j,{}_{-[i,m]}R^{j+1}\}\\&= \{{}_{-i}K:\
K\in\mathcal{K}_{i-1}{}^{\ast}\}- \{{}_{-[1,i]}R^{j-1},{}_{-[i,m]}R^{j+1}\}\ .
\end{split}
\end{equation*}
\end{itemize}

$\quad$

By induction, we have
\begin{equation*}
\mathcal{K}_s{}^{\ast}=\bmax^+\bigl(\mathfrak{V}(\pmb{R})\bigr)\ .
\end{equation*}
According to Proposition~{\rm\ref{p:4}}, this is a critical committee for $\mathcal{N}^s$.
\end{sketch}

\end{document}